\documentclass{amsart}[11 pt]

\usepackage{inputenc}
\usepackage{stmaryrd}
\input xy
\xyoption{all}
\usepackage{epsfig}
\usepackage{color}
\usepackage{amsthm}
\usepackage{changepage}
\usepackage{textcomp}
\usepackage{mathdots}
\usepackage{amsmath}
\usepackage{amsfonts}
\usepackage{amssymb}
\usepackage{mathrsfs}
\usepackage{amscd}
\usepackage{amsopn}
\usepackage{enumerate}
\usepackage{mathtools}
\usepackage{url}
\usepackage{graphicx}
\usepackage{tikz}
\usepackage{hyperref}
\usepackage{cite}
\usepackage{xcolor}
\usepackage[normalem]{ulem}
\usepackage{newunicodechar}
\usepackage{hyperref}\hypersetup{colorlinks}
\usepackage{color} 
\definecolor{darkred}{rgb}{1,0,0} 
\definecolor{darkgreen}{rgb}{0,0.8,0}
\definecolor{darkblue}{rgb}{0,0,1}
\hypersetup{colorlinks,
linkcolor=darkblue,
filecolor=darkgreen,
urlcolor=darkred,
citecolor=darkgreen}

\numberwithin{equation}{section}
\newtheorem {Theorem}{Theorem}
\numberwithin{Theorem}{section}

\newtheorem {Lemma}[Theorem] {Lemma}

\newtheorem {Proposition}[Theorem]{Proposition}
\newtheorem {Corollary}[Theorem]{Corollary}
\theoremstyle{definition}
\newtheorem{Definition}[Theorem]{Definition}

\newtheorem{Remark}[Theorem]{Remark}
\newtheorem{Example}[Theorem]{Example}

\theoremstyle{remark}

\expandafter\chardef\csname pre amssym.def at\endcsname=\the\catcode`\@
\catcode`\@=11
\def\undefine#1{\let#1\undefined}
\def\newsymbol#1#2#3#4#5{\let\next@\relax
\ifnum#2=\@ne\let\next@\msafam@\else
\ifnum#2=\tw@\let\next@\msbfam@\fi\fi
\mathchardef#1="#3\next@#4#5}
\def\mathhexbox@#1#2#3{\relax
\ifmmode\mathpalette{}{\m@th\mathchar"#1#2#3}%
\else\leavevmode\hbox{$\m@th\mathchar"#1#2#3$}\fi}
\def\hexnumber@#1{\ifcase#1 0\or 1\or 2\or 3\or 4\or 5\or 6\or 7\or 8\or
9\or A\or B\or C\or D\or E\or F\fi}
\font\teneufm=eufm10
\font\seveneufm=eufm7
\font\fiveeufm=eufm5
\newfam\eufmfam
\textfont\eufmfam=\teneufm
\scriptfont\eufmfam=\seveneufm
\scriptscriptfont\eufmfam=\fiveeufm

\catcode`\@=\csname pre amssym.def at\endcsname

\newcommand{\CE}{{\mathcal E}}

\newcommand{\fsl}{{\mathfrak {sl}}}

\newcommand{\sSL}{{\mathsf {SL}}}

\newcommand{\sG}{{\mathsf G}}

\newcommand{\sH}{{\mathsf H}}

\newcommand{\rB}{{\mathrm B}}

\newcommand{\rH}{{\mathrm H}}

\newcommand{\rU}{{\mathrm U}}

\newcommand{\rW}{{\mathrm W}}

\newcommand{\Hod}{{\mathrm {Hod}}}
\newcommand{\dR}{{\mathrm {dR}}}
\newcommand{\DeH}{{\mathrm {DH}}}
\newcommand{\VHS}{{\mathsf{VHS}}}
\newcommand{\NC}{{\mathsf{NC}}}

\newcommand{\NAH}{{\mathsf{NAH}}}

\newcommand{\Hom}{{\mathrm{Hom}}}

\newcommand{\End}{{\mathrm{End}}}
\newcommand{\Hit}{{\mathsf{Hit}}}

\newcommand{\oper}{{\mathsf{oper}}}

\newcommand{\ad}{{\mathrm{ad}}}

\newcommand{\Ll}{{\mathcal L}}

\newcommand{\Gg}{{\mathcal G}}
\newcommand{\Hh}{{\mathcal H}}

\newcommand{\Mm}{{\mathcal M}}

\newcommand{\ip}[2]{\langle #1, #2 \rangle}

\newcommand{\im}{\text{im}}
\newcommand{\dsp}{\displaystyle}
\newcommand{\lb}{\left\{}
\newcommand{\rb}{\right\}}

\def \C {{\mathbb C}}
\def \R {{\mathbb R}}
\def \Z {{\mathbb Z}}

\def \P {{\mathbb P}}

\def \lra {{\longrightarrow}}

\def \p {\partial}
\def \bp {\bar{\partial}}

\def \rk {\operatorname{rk}}

\def \im {\operatorname{im}}
\def \tr {\operatorname{tr}}

\def \CL {\operatorname{CL}}

\def \ind {\operatorname{ind}}

\newcommand{\An}{\xymatrix{ *{\circ} \ar@{-}[r]|*\dir{ } & *{\circ}\ar@{-}[r]&{\cdots}&*{\circ}\ar@{-}[l]|*\dir{ }\ar@{-}[r]|*\dir{ }&*{\circ} }}
\newcommand{\Anlabel}{\xymatrix@R=.25em{ *{\circ}\ar@<-1ex>@{}[d]^{\alpha_{1}} \ar@{-}[r]|*\dir{ } & *{\circ}\ar@<-1ex>@{}[d]^{\alpha_{2}}\ar@{-}[r]&{\cdots}&*{\circ}\ar@{-}[l]|*\dir{ }\ar@{-}[r]|*\dir{ }\ar@<-2ex>@{}[d]^{\alpha_{n-1}}&*{\circ}\ar@<-1ex>@{}[d]^{\alpha_{n}}\\&&&& }}
\newcommand{\Bn}{\xymatrix{ *{\circ} \ar@{-}[r]|*\dir{ } & *{\circ}\ar@{-}[r]&{\cdots}&*{\circ}\ar@{-}[l]|*\dir{ }\ar@{=}[r]|*\dir{>}&*{\circ} }}
\newcommand{\Bnlabel}{\xymatrix@R=.25em{ *{\circ}\ar@<-1ex>@{}[d]^{\alpha_{1}} \ar@{-}[r]|*\dir{ } & *{\circ}\ar@<-1ex>@{}[d]^{\alpha_{2}}\ar@{-}[r]&{\cdots}&*{\circ}\ar@{-}[l]|*\dir{ }\ar@{=}[r]|*\dir{>}\ar@<-2ex>@{}[d]^{\alpha_{n-1}}&*{\circ}\ar@<-1ex>@{}[d]^{\alpha_{n}}\\&&&& }}
\newcommand{\Cn}{\xymatrix{ *{\circ} \ar@{-}[r]|*\dir{ } & *{\circ}\ar@{-}[r]&{\cdots}&*{\circ}\ar@{-}[l]|*\dir{ }\ar@{=}[r]|*\dir{<}&*{\circ} }}
\newcommand{\Cnlabel}{\xymatrix@R=.25em{ *{\circ}\ar@<-1ex>@{}[d]^{\alpha_{1}} \ar@{-}[r]|*\dir{ } & *{\circ}\ar@<-1ex>@{}[d]^{\alpha_{2}}\ar@{-}[r]&{\cdots}&*{\circ}\ar@{-}[l]|*\dir{ }\ar@{=}[r]|*\dir{<}\ar@<-2ex>@{}[d]^{\alpha_{n-1}}&*{\circ}\ar@<-1ex>@{}[d]^{\alpha_{n}}\\&&&& }}
\newcommand{\Dn}{\xymatrix@R=.25em{&&&&*{\circ} \\ *{\circ} \ar@{-}[r]|*\dir{ } &
*{\circ}\ar@{-}[r]&{\cdots}&*{\circ}\ar@{-}[l]|*\dir{ }\ar@{-}[ur]|*\dir{ }\ar@{-}
[dr]|*\dir{ }& \\ &&&&*{\circ} }}
\newcommand{\Dnlabel}{\xymatrix@R=.25em{&&&&*{\circ}\ar@<-1ex>@{}[d]^{\alpha_{n-1}}
\\ *{\circ}\ar@<-1ex>@{}[d]^{\alpha_{1}} \ar@{-}[r]|*\dir{ } & *{\circ}\ar@<-1ex>@{}[d]^{\alpha_{2}}\ar@{-}[r]&{\cdots}&*{\circ}\ar@{-}[l]|*\dir{ }\ar@{-}[ur]|
*\dir{ }\ar@{-}[dr]|*\dir{ }\ar@<-4ex>@{}[d]^{\alpha_{n-2}}& \\ &&&&*{\circ}\ar@<-2ex>@{}[u]^(-.5){\alpha_{n}} }}

\usepackage[margin=1.1in]{geometry}

\title{The Upward Conformal limit}
\author[N. Alonso-Garcia]{Nerea Alonso-Garcia}
\address{DEPARTMENT OF MATHEMATICS, UNIVERSITY OF CALIFORNIA, RIVERSIDE, CA 92521, USA}
\email{nalon008@ucr.edu}

\begin{document}
\maketitle

\begin{abstract}
    For a compact Riemann surface $X$, both the Higgs bundle moduli space and the de Rham moduli space have stratifications by holomorphic Lagrangian fibrations over the space of complex variations of Hodge structure. The $\hbar$-conformal limit gives a biholomorphic identification of the fibers of the two strata over a stable complex variation of Hodge structure. We introduce here the $\hbar$-upward conformal limit, which is the limit at infinity rather than at zero. The main results of this paper relate the irreducible components of the nilpotent cone on $X$, which are Lagrangian submanifolds, with Lagrangians in the Morse strata on $\Bar{X}$, called zero velocity Lagrangians, via the $\hbar$-upward conformal limit. In particular, Higgs bundles that flow up to the zero of the Hitchin section on $X$ are identified with opers on $\Bar{X}$. 
\end{abstract}

\noindent

\section{Introduction}

For a compact Riemann surface $X$ of genus $g \geq 2$, the moduli space of Higgs bundles, $\Mm_\sH$, and the moduli space of holomorphic connections (or the de Rham moduli space), $\Mm_\dR$, are holomorphic symplectic manifolds. The first moduli space was introduced by Hitchin in \cite{MR887284} as solutions to the reduction to two dimensions of the self-dual Yang-Mills equations. Since then, it has been extensively studied and has become an important tool connecting different areas in mathematics and theoretical physics, such as symplectic geometry, representation theory, mirror symmetry \cite{MR1990670} and the geometric Langlands correspondence \cite{MR2306566}.

Recall that a Higgs bundle is a pair, $(E, \Phi)$, consisting of a holomorphic vector bundle $E$ on $X$, and an endomorphism valued holomorphic one form, $\Phi$, called the \textit{Higgs field}. There is a remarkable link between $\Mm_\sH$ and $\Mm_\dR$ which we describe next, where we identify $\Mm_\dR$ with the moduli space of flat connections on $X$. A choice of hermitian metric $h$ on $E$ associates hermitian adjoints to endormorphisms and ``completes'' the $\bp$-operator to the Chern connection: $\bp_E + \p_E^h$. When the Higgs bundle is polystable there is a special choice of metric, called the \textit{harmonic metric}, which produces a completely reducible flat connection. The map which builds this association is the celebrated nonabelian Hodge correspondence ($\NAH$) and it defines a homeomorphism between both moduli spaces:
\begin{eqnarray} \label{NAH}
    \begin{array}{cccl}
        \NAH \; :& \Mm_\sH & \xlongrightarrow{\sim} & \Mm_\dR \\
        &\left[(E, \Phi)\right] & \longmapsto & \left[ D = \bp_E + \p_E^h + \Phi + \Phi^{\ast_h}\right].
    \end{array}
\end{eqnarray}
This is a real analytic map, but it is not holomorphic, which means the moduli spaces are different as complex analytic manifolds. However, the complex structures combine via composition to define a hyperkähler structure. The corresponding twistor space is known as the Deligne-Hitchin moduli space and is denoted by $\Mm_\DeH$. From Deligne's perspective \cite{MR1492538}, $\Mm_\DeH$ is considered as the moduli space of $\lambda$\textit{-connections} on the Riemann surface $X$, $\Mm_\Hod$, glued with the moduli space of $\lambda$-connections on the conjugate Riemann surface $\Bar{X}$. These objects interpolate between Higgs bundles, when $\lambda = 0$, and flat connections, when $\lambda = 1$. More precisely, there is a holomorphic map $p: \Mm_{\Hod} \longrightarrow \C$ with two distinguished fibers: $p^{-1}(0) = \Mm_\sH$ and $p^{-1}(1) = \Mm_\dR$.

An important property of the twistor space is the existence of a $\C^\ast$-action extending the one defined by Hitchin on the Higgs bundle moduli space: $\xi \cdot [(\bp_E, \Phi)] = [(\bp_E, \xi \Phi)], \; \xi \in \C^\ast$. The Higgs bundles fixed by this action are the ones defining \textit{complex variations of Hodge structure} ($\VHS$). In \cite{MR2681730}, Simpson showed that the action defines a Białynicki-Birula type decomposition of $\Mm_{\Hod}$ where the strata $\rW_\alpha$ is labeled by the connected components $\VHS_{\alpha} \subset \Mm_{\sH}$ of the set of fixed points. This decomposition gives stratifications on the fibers $\Mm_\sH$ and $\Mm_\dR$. On the Higgs bundle side, the stratification, $\rW_\alpha^0 = \rW_\alpha \cap p^{-1}(0)$, is known as the \textit{Morse strata} defined by Hitchin; and in the de Rham side, $\rW_\alpha^1 = \rW_\alpha \cap p^{-1}(1)$, it is the \textit{partial oper strata} defined by Simpson. In both cases, each stratum is a holomorphic fibration by Lagrangians over a connected component of the set of complex variations of Hodge structure. Of special relevance are the Hitchin section and the set of opers which form a closed stratum of the stratifications on the spaces of Higgs bundles and holomorphic connections respectively. However, the nonabelian Hodge correspondence does not preserve these stratifications.

There is a different way of identifying the moduli space of Higgs bundles and the moduli space of holomorphic connections using conformal limits. For a stable Higgs bundle $(\bp_E, \Phi)$ with harmonic metric $h$, we can combine the $\C^\ast$-action on $\Mm_{\sH}$ and the \textit{real twistor line} in $\Mm_{\Hod}$ to scale the correspondence in (\ref{NAH}). More specifically,  for $R > 0$ the scaled Higgs bundle $(\bp_E, R \Phi)$ is also stable, and $\NAH$ produces harmonic metrics $h_R$. For $\xi \in \C^\ast$, the connection
\begin{eqnarray*}
    \bp_E + \p_E^{h_R} + \xi^{-1} R \Phi + \xi R \Phi^{\ast_{h_R}}
\end{eqnarray*}
is flat. Fixing the ratio $\hbar := \xi R^{-1}$, we obtain a one parameter family of flat connections:
\begin{eqnarray} \label{DR}
    D_{R, \hbar}(\bp_E, \Phi) = \bp_E + \p_E^{h_R} + \hbar^{-1} \Phi + \hbar R^2 \Phi^{\ast_{h_R}}.
\end{eqnarray}
When the limit $\displaystyle \lim_{R \to 0} D_{R, \hbar}(\bp_E, \Phi)$ exists, it is a flat connection called the $\hbar$-conformal limit.  The existence of this map for points in the Hitchin component was conjectured by Gaiotto in \cite{gaiotto2014operstba}, and later verified in \cite{MR4214341}. In \cite{MR3949609}, the existence of the $\hbar$-conformal limit was generalized to every Higgs bundle that flows down to a stable fixed point, i.e., Higgs bundles whose $\C^\ast$-limit at zero is stable. Furthermore, they showed that the conformal limit is defined on each stratum and holomorphically identifies the affine Lagrangian fibers of each stratum; in particular, it identifies Hitchin sections with the space of opers. Unlike the nonabelian Hodge correspondence, the conformal limit does not extend to a continuous map on the entire moduli space.

Another key feature of the Higgs moduli space is the Hitchin fibration. The \textit{nilpotent cone}, $\NC$, is the fiber over zero of this map, and has a complicated structure with several irreducible components. Each of these components, $\NC_{\alpha}$, is a Lagrangian submanifold \cite{MR887284, MR885778} and can be characterized as the set of points that flow up to a connected component of the set of complex variations of Hodge structure, $\VHS_{\alpha}$. The nilpotent cone has been an active object of study because it encodes important properties of the moduli space of Higgs bundles. In particular, we will show that it contains information about the Riemann surface with the opposite orientation.

The main results of this paper relate the irreducible components of the nilpotent cone on $X$ with the Morse stata on $\Bar{X}$. Instead of studying the limit as $R \to 0$ of the family of flat connections (\ref{DR}) for a stable Higgs bundle $(\bp_E, \Phi)$, we can look at the limit
\begin{eqnarray*}
    \dsp \lim_{R \to \infty} D_{R, \hbar}(\bp_E, \Phi),
\end{eqnarray*}
when it exists, it is a flat connection and we call it the \textit{$\hbar$-upward conformal limit}. Following a similar approach as in \cite{MR3949609}, we show that for Higgs bundles that flow up to stable fixed points, i.e., points in $\NC$ whose $\C^\ast$-limit at infinity is a stable fixed point, the $\hbar$-upward conformal limit exists. The requirement that the Higgs bundles being nilpotent is equivalent to the existence of their $\C^\ast$-limit at infinity.
\begin{Theorem} \label{MT}
    Let $[(\bp_E, \Phi)]$ be an element of the nilponent cone such that its limit at infinity under the $\C^\ast$-action is a stable Higgs bundle. Then the $\hbar$-upward conformal limit for $[(\bp_E, \Phi)]$ exists, and we denote it by $\CL_{\hbar,\infty}([(\bp_E, \Phi)])$.
\end{Theorem}

In general, each stratum $W^0_\alpha$ is not a Lagrangian submanifold, unless the corresponding $\VHS_\alpha$ consists of only one Higgs bundle. It is known that each $W^0_\alpha$ has a natural filtration by weights of the $\C^\ast$-action \cite[Section 4]{MR3949609}. In \cite{gothen2026zerovelocitylagrangianshitchin}, the authors use that filtration to identify a Lagrangian, $\Ll_\alpha$, in $W^0_\alpha$, which they call \textit{zero velocity Lagrangian}. A surprising consequence of the proof of Theorem \ref{MT} is that the $\hbar$-upward conformal limit of an irreducible component of the nilpotent cone in our Riemann surface $X$ is precisely the $\hbar^{-1}$-conformal limit of this zero velocity Lagrangian in $\Bar{X}$.
\begin{Theorem} \label{MT2}
    Let $\VHS_\alpha \subset \Mm_\sH^X$ be a connected component of the set of fixed points. Let $\NC_{\alpha}^X$ be the irreducible component of the nilpotent cone on $X$ and $\Ll_{\alpha}^{\Bar{X}}$ be the zero velocity Lagrangian submanifold of the Morse strata on $\Bar{X}$, both associated to $\VHS_\alpha$. Then,
    \begin{eqnarray*}
        \CL_{\hbar, \infty}^{X}(\NC_{\alpha}) = \CL_{\hbar^{-1}, 0}^{\Bar{X}}(\Ll_{\alpha}).
    \end{eqnarray*}
\end{Theorem}

\begin{Remark} \label{R 1}
    For any connected component $\VHS_\alpha$, there is a filtration of the Morse strata, $W_\alpha^0$, and the irreducible component of the nilpotent cone, $\NC_\alpha$, given by
    \begin{eqnarray*}
        W_\alpha^{0,\ell-1} \subset \dots \subset W_\alpha^{0,1} \subset W_\alpha^{0} \hspace{10pt} \mbox{ and } \hspace{10pt} \NC_\alpha^{1-\ell} \subset \dots \subset \NC_\alpha^{-2} \subset \NC_\alpha,
    \end{eqnarray*}
    where the components of the filtrations are labeled by the weight subspaces of the $\C^\ast$-action on $\Mm_\sH$. We prove the following general statement (see Theorem \ref{T 2}) which relates the above filtrations using conformal limits: for each $1 \leq j \leq \ell - 1$,
    \begin{eqnarray*}
        \CL_{\hbar, \infty}^{X}(\NC_{\alpha}^{-j}) = \CL_{\hbar^{-1}, 0}^{\Bar{X}}(W_\alpha^{0, j}).
    \end{eqnarray*}
    Observe that Theorem \ref{MT2} is the case $j = 1$.
\end{Remark}

As a consequence of the last result, we can identify the $\hbar$-upward conformal limit of Higgs bundles that flow up to a stable fixed point on $X$ with partial opers on $\Bar{X}$. In \cite{MR1174252}, Hitchin showed that for the unique fixed point in a Hitchin section, $[(\bp_\infty, \Phi_\infty)]$, the following holds
\begin{eqnarray*}
    W_\alpha^0(\bp_\infty, \Phi_\infty) \cong \bigoplus_{j = 1}^{n-1} H^0(K^{j+1}).
\end{eqnarray*}
Gaiotto conjectured \cite{gaiotto2014operstba} that the $\hbar$-conformal limit of points in $W_\alpha^0(\bp_\infty, \Phi_\infty)$ should be the space of opers. More precisely, for any $q \in \dsp \bigoplus_{j = 1}^{n-1} H^0(K^{j+1})$, the $\hbar$-conformal limit of the Higgs bundle associated to $q$ is the oper associated to $q$. This was proven in \cite{MR4214341}. In Section \ref{S 5}, we prove that for the unique fixed point in a Hitchin section, $[(\bp_\infty, \Phi_\infty)]$, we have
\begin{eqnarray*}
    \NC_\alpha(\bp_\infty, \Phi_\infty) \cong \bigoplus_{j = 1}^{n - 1} H^0(K_{\Bar{X}}^{j + 1}),
\end{eqnarray*}
which leads to the following result:
\begin{Corollary} \label{MT3}
    For any $\beta \in \dsp \bigoplus_{j = 1}^{n-1} H^1(K_{X}^{-j})$, the $\hbar$-upward conformal limit of the Higgs bundle associated to $\beta$ is the oper on $\Bar{X}$ associated to $\beta$ viewed as an element of $\dsp \bigoplus_{j = 1}^{n-1} H^0(K_{\Bar{X}}^{j + 1})$.
\end{Corollary}

\begin{Remark}
    One restriction made in the paper is that here we only consider Higgs vector bundles. However, analogous results should hold for $\sG$-Higgs bundles and flat $\sG$-connections for any semisimple Lie group $\sG$. Another restriction is that in all the main results we need the limiting $\VHS$ to be stable, it would be interesting to extend the statements of this paper to the polystable case using the techniques from \cite{kwong2025singularlagrangianshitchinmoduli}.
\end{Remark}

The paper is organized as follows: Section \ref{S 2} contains a brief review of the moduli spaces of Higgs bundles, flat connections and $\lambda$-connections, how they interact with each other, and some important features that are relevant for our purposes. In Section \ref{S 3} we construct a negative slice at a stable fixed point, and we use it to parametrize the unstable manifold at the fixed point in \ref{DF}. This machinery is needed to give a proof of the existence of the upward conformal limit in \ref{S 5}. In this last section, we also use the upward conformal limit to relate the Hitchin section and the nilpotent cone in the Riemann surface with opers and zero velocity Lagrangians in the Riemann surface with the opposite orientation.\\

\textit{Acknowledgements}. I would like to thank my advisor, Brian Collier, for suggesting this problem as well as for all the insightful conversations about background and ideas, the comments to improve the paper and his endless support. I would like to also thank André Oliveira and Rodrigo Pereira for helpful discussions. The author has been partially supported by NSF grant DMS-2337451.\\

\textit{AI statement}. No artificial intelligence tools were used in the research and writing of this paper.

\section{Background} \label{S 2}

Throughout the paper, we fix $X$ a compact Riemann surface of genus $g \geq 2$ and a smooth complex vector bundle $E \longrightarrow X$ of rank $n$ with a trivialization of the determinant line bundle $\det(E)$. In some parts of the text objects on the Riemann surface $X$ and its conjugate $\Bar{X}$ will appear, to distinguish between them in the notation we will indicate on which Riemann surface the object lives. For example, $K_X$ and $K_{\Bar{X}}$ denote the canonical line bundles on $X$ and $\Bar{X}$ respectively, but if we are only dealing with objects in $X$ then it will be denoted just by $K$.

\subsection{$\lambda$-connections and the Deligne-Hitchin moduli space}

We start by briefly reviewing two well-known examples of $\lambda$-connections: the moduli spaces of Higgs bundles and flat connections.

Recall that a holomorphic structure on the complex vector bundle $E \longrightarrow X$ is equivalent to a choice of $\bp$-operator, $\bp_E: \Omega(E) \longrightarrow \Omega^{0,1}(E)$, satisfying the Leibniz rule and inducing on $\det(E)$ the $\bp$-operator compatible with the trivialization.
\begin{Definition}
    An $\sSL(n, \C)$-Higgs bundle on $X$ is a pair $(\bp_E, \Phi)$ where $\Phi: E \longrightarrow E \otimes K$ is holomorphic ($\bp_E \Phi = 0$), and traceless ($\mbox{Tr}(\Phi) = 0$).
\end{Definition}
Let $\mathscr{H}$ denote the set of all $\sSL(n, \C)$-Higgs bundles on $E$.
A Higgs bundle $(\bp_E, \Phi)$ is called \textit{stable} (resp. \textit{semistable}) if for all $\Phi$-invariant holomorphic subbundles $F \subset (E, \bp_E)$ with $0 < \rk(F) < \rk(E)$, we have $\deg(F) < 0$ (resp. $\deg(F) \leq 0$). The Higgs bundle is \textit{polystable} if it is semistable and a direct sum of Higgs bundles $(\bp_{E_j}, \Phi_j)$ so that $\deg(E_j) = 0$ for all $j$ and each $(\bp_{E_j}, \Phi_j)$ is stable. The gauge group $\Gg(E)$ of smooth automorphisms of $E$ that act trivially on $\det(E)$ acts on the set of $\sSL(n,\C)$-Higgs bundles by the adjoint action: $(\bp_E, \Phi) \cdot g = (g^{-1} \circ \bp_E \circ g, g^{-1} \circ \Phi \circ g)$.  We define the \textit{Higgs bundle moduli space} by
\begin{eqnarray*}
    \Mm_\sH = \mathscr{H}^{ps} / \Gg(E),
\end{eqnarray*}
where $\mathscr{H}^{ps}$ is the set of polystable Higgs bundles. Its smooth locus consists of equivalence classes of stable Higgs bundles, and $\Mm_\sH$ is a (singular) complex analytic space which is a normal, quasiprojective variety of dimension $(n^2 -1)(2g-2)$ \cite{MR1085642, MR887284, MR1320603}.\\

The gauge group also acts on $\mathscr{D}$, the space of flat $\sSL(n, \C)$-connections that induce the trivial connection on $\det(E)$, by $D \cdot g = g^{-1} \circ D \circ g$. A connection $D$ is called \textit{completely reducible} if every $D$-invariant subbundle $F \subset E$ has a $D$-invariant complement, and $D$ is called \textit{irreducible} if there are no nontrivial $D$-invariant subbundles. We define the \textit{de Rham moduli space} by
\begin{eqnarray*}
    \Mm_\dR = \mathscr{D}^{cr}/\Gg(E),
\end{eqnarray*}
where $\mathscr{D}^{cr}$ denotes the set of completely reducible flat connections. The smooth locus of $\Mm_\dR$ consists of equivalence classes of irreducible representations, and $\Mm_\dR$ is a (singular) complex analytic space which is Stein and has complex dimension $(n^2 -1)(2g-2)$ \cite{MR1320603}.\\

The famous nonabelian Hodge correspondence \cite{MR965220, MR887285, MR887284, MR944577} gives a relation between Higgs bundles and flat connections. We mentioned that both moduli spaces are complex manifolds; however, their complex structures are not identified via the nonabelian Hodge correspondence. Instead, this map provides $\Mm_\dR$ with a hyperkähler structure. In order to give an analytic description of its twistor space, $\Mm_\dR \times \P^1$, Deligne introduced the concept of a $\lambda$-connection. These objects have been studied by Simpson \cite{MR1492538}, and many others.

\begin{Definition}
    Let $\lambda \in \C$. A \textit{holomorphic $\lambda$-connection} (or just \textit{$\lambda$-connection}) on $E$ is a triple $(\lambda, \bp_E, \nabla_\lambda)$ where
    \begin{itemize}
        \item $\bp_E$ is a $\bp$-operator on $E$,

        \item $\nabla_\lambda: \Omega(E) \longrightarrow \Omega^{1,0}(E)$ is a $C^\infty$ differential operator satisfying
        \begin{eqnarray*}
            \nabla_\lambda(fs) = \lambda \p f \otimes s + f \cdot \nabla_\lambda s, \; \mbox{ for all } f \in C^\infty(X), s \in \Omega^0(E),
        \end{eqnarray*}
        such that the induced differential operator on $\det(E)$ coincides with $\lambda \p_E$.

        \item $\nabla_\lambda$ is compatible with the holomorphic structure, i.e., $[\bp_E, \nabla_\lambda] = 0$.
        \end{itemize}
\end{Definition}

Observe that $\lambda$-connections give a way of interpolating between flat connections and Higgs bundles. Indeed, when $\lambda = 1$, a 1-connection is just a holomorphic connection $\nabla : \CE \longrightarrow \CE \otimes K$ on the holomorphic bundle $\CE = (E, \bp_E)$ defined by $\bp_E$, which is equivalent to a flat connection $D = \bp_E + \nabla$ on $E$. In general, for $\lambda \neq 0$ a $\lambda$-connection $(\lambda, \bp_E, \nabla_\lambda)$ gives us a flat $\sSL(n,\C)$-connection on $E$ via $D = \bp_E + \lambda^{-1} \nabla_\lambda$. When $\lambda = 0$, then $\Phi:=\nabla_0: \CE \longrightarrow \CE \otimes K$ is linear and holomorphic, thus $(\bp_E, \Phi)$ defines an $\sSL(n, \C)$-Higgs bundle.

The gauge group $\Gg(E)$ acts on the space of $\lambda$-connections by
\begin{eqnarray*}
    (\lambda, \bp_E, \nabla_\lambda) \cdot g = (\lambda, g^{-1} \circ \bp_E \circ g, g^{-1} \circ \nabla_\lambda \circ g).
\end{eqnarray*}
Note that this action specializes to the usual actions on the space of Higgs bundles and flat connections that we described above.

To define the moduli space of $\lambda$-connections we need the notion of stability.

\begin{Definition}
    A holomorphic $\lambda$-connection on $E$ is called
    \begin{itemize}
        \item \textit{semistable} if for all $\nabla_\lambda$-invariant holomorphic subbundles $F \subset (E, \bp_E)$ we have $\deg F \leq 0$,
            
        \item \textit{stable} if for all $\nabla_\lambda$-invariant holomorphic subbundles $F \subset (E, \bp_E), 0 < \rk(F) \leq \rk(E)$, we have $\deg F <0$,

        \item \textit{polystable} if it is semistable and a direct sum of of stable $\lambda$-connections whose associated holomorphic bundles have degree zero.
    \end{itemize}
\end{Definition}

\begin{Remark}
    When $\lambda = 0$, the above definition recovers the notion of a stable (resp. semistable, resp. polystable) Higgs bundle. Similarly, when $\lambda = 1$, a 1-connection is stable (resp. polystable) if and only if the associated flat connection is irreducible (resp. completely reducible).
\end{Remark}

The moduli space of (polystable) $\lambda$-connections or \textit{Hodge moduli space} is defined by
\begin{eqnarray*}
    \Mm_\Hod(X, \sSL(n, \C)) := \{(\lambda, \bp_E, \nabla_\lambda) \mbox{ polystable holomorphic } \lambda\mbox{-connection} \}\ / \Gg(E).
\end{eqnarray*}
In general, we just denote this moduli space by $\Mm_\Hod$, and we will specify the dependence on the Riemann surface when various Riemann surfaces are involved in the discussion.
This is a complex space, and its smooth locus is precisely the locus of stable $\lambda$-connections.

There is a $\C^\ast$-action on $\Mm_\Hod$ given by
\begin{eqnarray} \label{action}
    \xi \cdot [(\lambda, \bp_E, \nabla_{\lambda})] = [(\xi \lambda, \bp_E, \xi \nabla_{\lambda})].
\end{eqnarray}
The combination of this action with the nonabelian Hodge map from (\ref{NAH}) produces a family of flat connections parametrized by $\xi \in \C$:
\begin{eqnarray*}
    \bp_E + \p_E^h + \xi^{-1} \Phi + \xi \Phi^{\ast_h}.
\end{eqnarray*}
Using the $\C^\ast$-action on $\lambda$-connections we get
\begin{eqnarray} \label{Conec}
    \xi \cdot [(1, \bp_E + \xi \Phi^{\ast_h}, \p_E^h + \xi^{-1} \Phi)] = [(\xi, \bp_E + \xi \Phi^{\ast_h}, \xi \p_E^h + \Phi)].
\end{eqnarray}
Moreover,
\begin{eqnarray} \label{limit}
    \lim_{\xi \to 0}[(\xi, \bp_E + \xi \Phi^{\ast_h}, \xi \p_E^h + \Phi)] = [(\bp_E, \Phi)],
\end{eqnarray}
which allows us to extend the $\C^\ast$-family (\ref{Conec}) to $\C$, this is an example of a \textit{real twistor line} in $\Mm_\Hod$. Thus, the above real twistor line through a Higgs bundle $[(\bp_E, \Phi)]$ is a section of $\Mm_\Hod$ which interpolates between the Higgs bundle and the flat connection given by $\NAH(\bp_E, \Phi)$.

\begin{Definition} \label{DeH}
    The \textit{Deligne-Hitchin moduli space} is
    \begin{eqnarray*}
        \Mm_{\DeH} = \Mm_\Hod^X \; \dot{\cup} \; \Mm_\Hod^{\Bar{X}}/\sim_{\DeH},
    \end{eqnarray*}
    where $(\lambda, \bp, \nabla) \sim_{\DeH} (\lambda^{-1}, \lambda^{-1}\bp, \lambda^{-1}\nabla)$ for any $(\lambda, \bp, \nabla) \in \Mm_\Hod^X$ with $\lambda \neq 0$.
\end{Definition}
The Deligne-Hitchin moduli space gives an analytic description of the twistor space $\Mm_\dR \times \P^1$ \cite{MR1492538}, but we will not need this perspective. For more details, see \cite{MR3927087, MR4297184}.

\subsection{Fixed points of the $\C^\ast$-action}

The holomorphic projection $p: \Mm_\Hod \longrightarrow \C$ given by
\begin{eqnarray*}
    p([(\lambda, \bp_E, \nabla_\lambda)]) = \lambda,
\end{eqnarray*}
has two special fibers: $p^{-1}(0) = \Mm_\sH$ and $p^{-1}(1) = \Mm_\dR$. The action described in (\ref{action}) only preserves the fiber $p^{-1}(0)$, and the restriction to this fiber is the $\C^\ast$-action on Higgs bundles:
\begin{eqnarray*}
    \xi \cdot [(\bp_E, \Phi)] = [(\bp_E, \xi \Phi)],
\end{eqnarray*}
i.e., the group $\C^\ast$ acts holomorphically on $\Mm_\sH$ by scaling the Higgs field. In particular, $\dsp \lim_{\xi \to 0} [(\bp_E, \xi \Phi)]$ always exists by properness of the Hitchin fibration (see Section \ref{S 2.3}) and is a fixed point.

\begin{Definition}
    A Higgs bundle $(\bp_E, \Phi)$ is called a \textit{complex variation of Hodge structure} ($\VHS$) if $[(\bp_E, \Phi)]$ is a fixed point of the $\C^\ast$-action on $\Mm_\sH$. We will denote by $\VHS_\alpha$ a connected component of the set of fixed points, labeled by some index $\alpha$.
\end{Definition}

Let $\dsp \amalg_\alpha \VHS_\alpha$ be the decomposition of $\VHS$ into connected components. In \cite{MR1492538}, Simpson defined a Białynicki-Birula type stratification of $\Mm_\Hod$ using the limit from (\ref{limit}). Let
\begin{eqnarray*}
    W_\alpha := \{ [(\lambda, \bp_E, \nabla_\lambda)] \in \Mm_\Hod : \lim_{\xi \to 0} [(\xi \lambda, \bp_E, \xi \nabla_\lambda)] \in \VHS_\alpha \}.
\end{eqnarray*}
Then $\Mm_\Hod$ decomposes as $\amalg_\alpha W_\alpha$, where each $W_\alpha$ is connected and foliated by $\C^\ast$-orbits. The restrictions $W^0_\alpha = W_\alpha \cap p^{-1}(0)$ and $W^1_\alpha = W_\alpha \cap p^{-1}(1)$ define stratifications of $\Mm_\sH$ and $\Mm_\dR$ respectively. However, these stratifications are not preserved by the nonabelian Hodge correspondence.\\

The following result gives a characterization of the fixed points \cite[Lemma 4.1]{MR1179076}.

\begin{Proposition} \label{P FP}
    A point $[(\bp_E, \Phi)] \in \Mm_\sH$ is a fixed point of the $\C^\ast$-action if and only if there is a smooth splitting $E = E_1 \oplus \dots \oplus E_\ell$ with respect to which
    \begin{eqnarray} \label{E FP}
        \bp_E = \begin{pmatrix}
        \bp_{E_1} &  & \\
          & \ddots & &  \\
          &         & &\bp_{E_\ell}
        \end{pmatrix} \hspace{10pt} \mbox{ and } \hspace{10pt} \Phi = \begin{pmatrix}
        0 & & &\\
         \Phi_1 & 0 & &\\
          & \ddots & \ddots & \\
         &         & \Phi_{\ell - 1} &0
    \end{pmatrix},
    \end{eqnarray}
    where $\Phi_j: E_j \longrightarrow E_{j+1} \otimes K$ is holomorphic. Observe that when $\ell = 1$, then $\Phi = 0$.
\end{Proposition}

The next result due to Simpson \cite{MR1179076} characterizes the flat connections associated to polystable Higgs bundles which are also fixed points.

\begin{Proposition} \label{Simpson}
    If $(\bp_E, \Phi)$ is a polystable Higgs bundle of the form (\ref{E FP}), then the splitting is orthogonal with respect to the harmonic metric $h$. In particular, the associated flat connection $D = \bp_E + \p_E^h + \Phi + \Phi^{\ast_h}$ is given by
    \begin{eqnarray*}
        \begin{pmatrix}
        \bp_{E_1} & \Phi_1^{\ast_h} &  & \\
          & \ddots & \ddots & &  \\
          & & \bp_{E_{\ell -1}} & \Phi_{\ell - 1}^{\ast_h}\\
          &         & &\bp_{E_\ell}
        \end{pmatrix} + \begin{pmatrix}
        \p_{E_1}^h & & &\\
         \Phi_1 & \p_{E_2}^h & &\\
          & \ddots & \ddots & \\
         &         & \Phi_{\ell - 1} & \p_{E_\ell}^h
    \end{pmatrix},
    \end{eqnarray*}
    where $\Phi_j^{\ast_h}: E_{j+1} \longrightarrow E_j \otimes \Bar{K}$ is the adjoint of $\Phi_j$ with respect to the induced metrics on $E_j$ and $E_{j+1}$.
\end{Proposition}

\begin{Example}
    An example of general fixed points are the so called \textit{(principal) Fuchsian} or \textit{uniformizing} Higgs bundles:
    \begin{eqnarray} \label{Fuchsian}
        (E, \bp_0, \Phi_0) = \left( L \oplus (L \otimes K^{-1}) \oplus \dots (L \otimes K^{n - 1}), \begin{pmatrix}
            0 &  & &\\
            1 & 0 & &\\
            & \ddots & \ddots &\\
            & & 1 & 0
        \end{pmatrix}\right),
    \end{eqnarray}
    where $L = K^{\frac{n - 1}{2}}$ is a choice of line bundle. When $n = 2$, the uniformizing Higgs bundle is
    \begin{eqnarray} \label{Fuchsian2}
        \left( K^{\frac{1}{2}} \oplus K^{-\frac{1}{2}}, \begin{pmatrix}
            0 & 0\\
            1 & 0
        \end{pmatrix}\right),
    \end{eqnarray}
    whose associated flat connection corresponds to the Fuchsian representation that uniformizes the Riemann surface $X$ \cite{MR887284}. Taking the $(n-1)$st-symmetric product, renormalizing and twisting by an $n$th-root of the trivial bundle one obtains (\ref{Fuchsian}). Observe that they are all stable Higgs bundles and, by Proposition \ref{Simpson}, their associated flat connection is given by
    \begin{eqnarray*}
        \begin{pmatrix}
            \bp_L & 1^{\ast_h} & & & \\
            & \ddots & \ddots & & \\
            & & \bp_{L \otimes K^{n - 2}} & 1^{\ast_h} &\\
            & & & \bp_{L \otimes K^{n - 1}}
        \end{pmatrix} + 
        \begin{pmatrix}
            \p_L^h & & & & \\
            1 & \p_{L \otimes K}^h&& & \\
            & \ddots & \ddots & &\\
            & & 1 & \p_{L \otimes K^{n - 1}}^h
        \end{pmatrix}.
    \end{eqnarray*}
\end{Example}

The splitting of $E$ in Proposition \ref{P FP} gives a $\Z$-grading on $\End(E)$ (for us traceless endomorphisms since we are working with $\sSL(n, \C)$-Higgs bundles): 
\begin{eqnarray} \label{Grading}
    \End(E) = \dsp \bigoplus_{j \in \Z} \End_j(E) \hspace{3pt} \mbox{ with } \hspace{3pt} \End_j(E) := \left\{ \begin{array}{ll}
       \dsp \bigoplus_{i - k = j}\Hom(E_i, E_k),  & \mbox{if } 1 - \ell \leq j \leq \ell - 1  \\
       \{0\} & \mbox{otherwise}
    \end{array}.
    \right.
\end{eqnarray}
We will use the same notation as in \cite{MR3949609}:
\begin{eqnarray} \label{Grading 2}
    N_+ := \bigoplus_{j > 0} \End_j(E), \hspace{11pt} N_- := \bigoplus_{j < 0} \End_j(E), \hspace{10pt} \mbox{ and } \hspace{10pt} L := \End_0(E) \cap \fsl(E).
\end{eqnarray}

The following lemma is immediate and will be useful later:
\begin{Lemma} \label{L 4}
    Let $(\bp_E, \Phi) \in \VHS_\alpha^X$ be stable with harmonic metric $h$. Then, $(\p_E^h, \Phi^{\ast_h}) \in \VHS_\alpha^{\Bar{X}}$ is stable with harmonic metric $h$. Moreover, $\End_j^{\Bar{X}}(E) = \End_{-j}^X(E)$ for any $j \in \Z$, and
    \begin{eqnarray*}
        N_+^{\Bar{X}} = N_-^X, \hspace{10pt} N_-^{\Bar{X}} = N_+^X, \hspace{9pt} \mbox{ and } \hspace{9pt} L^{\Bar{X}} = L^X.
    \end{eqnarray*}
\end{Lemma}

\subsection{Hitchin section and opers} \label{S 2.3}

An important feature of the Higgs bundle moduli space is the Hitchin fibration, which esentially maps a Higgs bundle $[(\bp_E, \Phi)]$ into the coefficients of its characteristic polynomial $\Phi$.
This map is proper and it turns $\Mm_\sH$ into a complete integrable system \cite{MR1174252, MR885778}. The preimage of zero of the Hitchin fibration is called the \textit{nilpotent cone}, $\NC$, because $[(\bp_E, \Phi)] \in \NC$ if and only if $\Phi$ is nilpotent. The connected components $\VHS_\alpha$ of the set of fixed points give a decomposition of the nilpotent cone in terms of its irreducible components: $\NC = \amalg_\alpha \NC_\alpha$, with $\NC_\alpha = W^\infty_\alpha$ \cite{MR1650276, MR4519161} where
\begin{eqnarray*}
    W^\infty_\alpha := \{ [(\lambda, \bp_E, \nabla_\lambda)] \in \Mm_\Hod : \lim_{\xi \to \infty} [(\xi \lambda, \bp_E, \xi \nabla_\lambda)] \in \VHS_\alpha \}.
\end{eqnarray*}
Each $\NC_\alpha$ is a Lagrangian subvariety \cite{MR887284, MR885778}.

The property of the fibration relevant to us is the existence of a section, called the \textit{Hitchin section}, 
\begin{eqnarray*}
    \begin{array}{rcl}
       p_{\Hit}: \dsp \bigoplus_{j=1}^{n-1} H^0(K^{j+1}) & \longrightarrow & \Mm_\sH \\
       (q_2, \dots, q_n) & \longmapsto & [(\bp_0, \Phi_0 + q_2 \otimes e_1 + \dots + q_n \otimes e_{n-1})] 
    \end{array},
\end{eqnarray*}
where $[(\bp_0, \Phi_0)] \in \VHS$ is the Higgs bundle described in (\ref{Fuchsian}), and $e_j$ is the $n \times n$ matrix with ones in the $j$th-superdiagonal. For any choice of line bundle $L$ the image of this map is called a \textit{Hitchin component}, and it defines a holomorphic Lagrangian submanifold on $\Mm_\sH$ \cite{MR1174252}.\\

There is also a distinguished holomorphic Lagrangian submanifold on the de Rham side, the set of \textit{opers}. It is well-known that opers generalize complex projective structures on $X$, and they have relevant applications to the geometric Langlands' program \cite{MR760998, beilinson2005opers}.
\begin{Definition}
    An $\sSL(n, \C)$-oper on $X$ is a triple $(E, \nabla, F_\bullet)$ where
    \begin{itemize}
        \item $E$ is a smooth complex vector bundle of rank $n$ over $X$, equipped with a trivialization of the determinant bundle $\det(E)$,

        \item $\nabla$ is a flat connection on $E$,

        \item $F_\bullet$ is a filtration of $E$ by subbundles, $0 = F_0 \subset F_1 \subset \dots \subset F_N = E$,
    \end{itemize}
    such that
    \begin{itemize}
        \item each $F_j$ is holomorphic (with respect to the holomorphic structure $\bp_\nabla$),

        \item $\nabla(F_j) \subset F_{j + 1} \otimes K$ for all $0 \leq j \leq N -1$,

        \item the induced linear map
        \begin{eqnarray*}
            \nabla : F_j/F_{j - 1} \longrightarrow (F_{j + 1}/F_j) \otimes K,
        \end{eqnarray*}
        is an isomorphism of line bundles for $1 \leq j \leq N - 1$. 
    \end{itemize}
\end{Definition}
Similarly to the Hitchin section, the space of opers is parametrized by the vector space $\dsp \bigoplus_{j=1}^{n-1} H^0(K^{j+1})$:
\begin{eqnarray} \label{oper}
    \begin{array}{rcl}
       p_{\oper}: \dsp \bigoplus_{j=1}^{n-1} H^0(K^{j+1}) & \longrightarrow & \Mm_\dR \\
       (q_2, \dots, q_n) & \longmapsto & [(\bp_0 + \Phi_0^{\ast_h}, \p_0^h + \Phi_0 + q_2 \otimes e_1 + \dots + q_n \otimes e_{n - 1})]
    \end{array},
\end{eqnarray}
where $[(\bp_0, \Phi_0)]$ and $e_j$ are as above.\\

Note that even when $n = 2$, nonabelian Hodge sends points in the Hitchin section to connections with holonomy on $\sSL(n, \R)$, but the holonomy of an oper is not real in general (see \cite[Example 3.8]{MR3949609}). This means that nonabelian Hodge does not map the Hitchin section into the set opers. However, in \cite{MR4214341}, the authors proved that there exists a biholomorphism between these two Lagrangians given by the $\hbar$-conformal limit.

\section{The negative slice at a $\VHS$} \label{S 3}

The purpose of this section is to define and understand the negative slice at a stable fixed point. This object will be relevant in Section \ref{DF}, where we will use it to parametrize Higgs bundles that flow up to a given fixed point. First, we recall the local Kuranishi model around a stable Higgs bundle $[(\bp_E, \Phi)] \in \Mm_\sH$, for more details on this and a general approach see \cite[Chapter 7]{MR909698}. We denote by $h$ the harmonic metric for $(\bp_E, \Phi)$. Throughout the rest of the paper, an unannotated $\ast$ will refer to adjoints with respect to $h$, and we set $\p_E^h := \p_E$. Then, the connection $D = \bp_E + \p_E + \Phi + \Phi^\ast$ given by the nonabelian Hodge correspondence is flat. Set
\begin{eqnarray*}
    D:= D'' + D', \hspace{1cm} D'' := \bp_E + \Phi \hspace{1cm}  \mbox{ and } \hspace{1cm} D' := \p_E + \Phi^\ast.
\end{eqnarray*}

Recall that
\begin{eqnarray*}
    T_{[(\bp_E, \Phi)]} \Mm_\sH = \frac{T_{(\bp_E, \Phi)} \mathscr{H}^{ps}}{T_{(\bp_E, \Phi)} \Gg(\bp_E, \Phi)} ,
\end{eqnarray*}
where $\Gg(\bp_E, \Phi)$ is the $\Gg(E)$-orbit of $(\bp_E, \Phi)$. We can realize both tangent spaces in the quotient above, $T_{(\bp_E, \Phi)} \mathscr{H}^{ps}$ and $T_{(\bp_E, \Phi)} \Gg(\bp_E, \Phi)$, as the kernel and image of $D''$ respectively. In particular, the deformation complex of $(\bp_E, \Phi)$ is given by
\begin{eqnarray} \label{DC}
    C(\bp_E, \Phi) : \Omega^0(\fsl(E)) \xlongrightarrow{D''} \Omega^{0,1}(\fsl(E)) \oplus \Omega^{1,0}(\fsl(E)) \xlongrightarrow{D''} \Omega^{1,1}(\fsl(E)).
\end{eqnarray}
Therefore
\begin{eqnarray*}
    T_{[(\bp_E, \Phi)]} \Mm_\sH = \frac{\ker (D'': \Omega^{0,1}(\fsl(E)) \oplus \Omega^{1,0}(\fsl(E)) \longrightarrow \Omega^{1,1}(\fsl(E)))}{\im (D'': \Omega^0(\fsl(E)) \longrightarrow \Omega^{0,1}(\fsl(E)) \oplus \Omega^{1,0}(\fsl(E)))} = H^1(C(\bp_E, \Phi)).
\end{eqnarray*}
Moreover, the normal space to $T_{(\bp_E, \Phi)} \Gg(\bp_E, \Phi)$ is
\begin{eqnarray*}
    \{(\beta, \varphi) \in \Omega^{0,1}(\fsl(E)) \oplus \Omega^{1,0}(\fsl(E)): (D'')^\ast (\beta, \varphi) = 0\}.
\end{eqnarray*}
Thus, $H^1(C(\bp_E, \Phi))$ is isomorphic to
\begin{eqnarray*}
     \Hh^1(\bp_E, \Phi) &:=& \lb (\beta, \varphi) \in \Omega^{0,1}(\fsl(E)) \oplus \Omega^{1,0}(\fsl(E)): D''(\beta, \varphi) = 0 \mbox{ and } (D'')^\ast (\beta, \varphi) = 0 \rb\\
     &=& \lb (\beta, \varphi) \in \Omega^{0,1}(\fsl(E)) \oplus \Omega^{1,0}(\fsl(E)): D''(\beta, \varphi) = 0 \mbox{ and } D' (\beta, \varphi) = 0 \rb,
\end{eqnarray*}
since $\ker (D'')^\ast = \ker D'$ because the deformation complex (\ref{DC}) is elliptic, and $D''$ and $D'$ satisfy the Kähler identities, in particular $(D'')^\ast = -i[\Lambda, D']$.
Hence
\begin{eqnarray*}
    T_{[(\bp_E, \Phi)]} \Mm_\sH \cong \Hh^1(\bp_E, \Phi).
\end{eqnarray*}

\begin{Definition}
    The \textit{Hodge slice} (or simply \textit{slice}) at $(\bp_E, \Phi) \in \mathscr{H}$ is given by
    \begin{eqnarray*}
        S_{(\bp_E, \Phi)} = \lb (\beta, \varphi) \in \Omega^{0,1}(\fsl(E)) \oplus \Omega^{0,1}(\fsl(E)) : D''(\beta, \varphi) + [\beta, \varphi] = 0 \mbox{ and } D'(\beta, \varphi) = 0 \rb.
    \end{eqnarray*}
\end{Definition}
Note that the first equation in the definition of the Hodge slice implies that $(\bp_E + \beta, \Phi + \varphi)$ is a Higgs bundle, so $S_{(\bp_E, \Phi)} \subset \mathscr{H}$. Recall that if $(\bp_E, \Phi)$ is stable, the slice defines a local chart in $\Mm_\sH$ through the projection map \cite[Theorem 7.3.17]{MR909698}:
\begin{eqnarray} \label{proj1}
    \begin{array}{rccl}
        p_\rH :& S_{(\bp_E, \Phi)} \cap \mathscr{H}^s & \longrightarrow & \Mm_\sH \\
        &(\beta, \varphi) & \longmapsto & [(\bp_E + \beta, \Phi + \varphi)].
    \end{array}
\end{eqnarray}
Observe that for any $\lambda \in \C$, we can extend the above map to the moduli space of $\lambda$-connections:
\begin{eqnarray*}
    \begin{array}{rccl}
       p_\lambda :& S_{(\bp_E, \Phi)} & \longrightarrow & \Mm_\Hod \\
        &(\beta, \varphi) & \longmapsto & [(\lambda, \bp_E + \lambda \Phi^\ast + \beta, \lambda \p_E + \Phi + \varphi)].
    \end{array}
\end{eqnarray*}

For the rest of the section we assume that $(\bp_E, \Phi)$ is a $\VHS$ and recall the grading associated to it from (\ref{Grading}). This gives a decomposition of the tangent space into weight subspaces, and we can decompose the deformation complex (\ref{DC}) into its graded subcomplexes:
\begin{eqnarray*}
    C_j(\bp_E, \Phi) : \Omega^0(\End_j(E)) \xlongrightarrow{D''} \Omega^{0,1}(\End_j(E)) \oplus \Omega^{1,0}(\End_{j-1}(E)) \xlongrightarrow{D''} \Omega^{1,1}(\End_{j-1}(E)),
\end{eqnarray*}
since $\bp_E$ preserves the grading and $\ad_{\Phi}: \End_j(E) \longrightarrow \End_{j-1} \otimes K$ because $\Phi \in H^0(\End_{-1} \otimes K)$. This also gives a grading $\Hh^1(\bp_E, \Phi) = \dsp \bigoplus_{j= 1 - \ell}^{\ell} \Hh^1_j(\bp_E, \Phi)$, where $\Hh^1_j(\bp_E, \Phi) := \ker D'' \cap \ker D'$ for the middle term in $C_j(\bp_E, \Phi)$. Set
\begin{eqnarray*}
    \Hh^1_+(\bp_E, \Phi) := \dsp \bigoplus_{j=1}^{\ell} \Hh^1_j(\bp_E, \Phi), \hspace{1cm} \mbox{ and } \hspace{1cm} \Hh^1_-(\bp_E, \Phi) := \dsp \bigoplus_{j=1 - \ell}^{-1} \Hh^1_j(\bp_E, \Phi).
\end{eqnarray*}
With the notation of (\ref{Grading 2}), $\Hh^i_+(\bp_E, \Phi) \cong H^i(C_+(\bp_E, \Phi))$ and $\Hh^i_-(\bp_E, \Phi) \simeq H^i(C_-(\bp_E, \Phi))$ for the subcomplexes of (\ref{DC}) given by
\begin{eqnarray*}
    C_+(\bp_E, \Phi) : \Omega^0(N_+) \xlongrightarrow{D''} \Omega^{0,1}(N_+) \oplus \Omega^{1,0}(L \oplus N_+) \xlongrightarrow{D''} \Omega^{1,1}(L \oplus N_+),
\end{eqnarray*}
and
\begin{eqnarray*}
    C_-(\bp_E, \Phi): \Omega^0(N_-) \xlongrightarrow{D''} \Omega^{0,1}(N_-) \oplus \Omega^{1,0}(N_{<-1}) \xlongrightarrow{D''} \Omega^{1,1}(N_{<-1}),
\end{eqnarray*}
respectively.

In \cite{MR3949609}, the authors proved that the $\rB \rB$-slice at $(\bp_E, \Phi)$,
\begin{eqnarray*}
    S^+_{(\bp_E, \Phi)} := \{(\beta, \varphi) \in \Omega^{0,1}(N_+) \oplus \Omega^{1,0}(L \oplus N_+) : D''(\beta, \varphi) + [\beta, \varphi] = 0, \, D'(\beta, \varphi) = 0\},
\end{eqnarray*}
defines a global coordinate chart through the Kuranishi map:
\begin{eqnarray*}
    \begin{array}{rccl}
        k:& \Omega^{0,1}(\fsl(E)) \oplus \Omega^{1,0}(\fsl(E)) & \xlongrightarrow{\sim} & \Omega^{1,1}(\fsl(E)) \\
        &(\beta, \varphi) & \longmapsto & k(\beta, \varphi) := (\beta, \varphi) + (D'')^\ast G([\beta, \varphi]),
    \end{array}
\end{eqnarray*}
where $G$ is the Green's operator for the Laplacian $D''(D'')^\ast$ on $\Omega^{1,1}(\fsl(E))$. Observe that when $(\beta, \varphi)$ is in $ S^+_{(\bp_E, \Phi)}$ then $k(\beta, \varphi)$ is an element of $\Hh^1_+(\bp_E, \Phi)$. 
We will show an analogous result for the negative slice and $\Hh^1_-(\bp_E, \Phi)$ following the same argument as in \cite{MR3949609}.
\begin{Definition}
    The negative slice at $(\bp_E, \Phi)$ is defined as
    \begin{eqnarray*}
        S^{<0}_{(\bp_E, \Phi)} = \{(\beta, \varphi) \in \Omega^{0,1}(N_-) \oplus \Omega^{1,0}(N_{<-1}) : D''(\beta, \varphi) + [\beta, \varphi] = 0, \, D'(\beta, \varphi) = 0\}.
    \end{eqnarray*}
\end{Definition}

\begin{Lemma}
    The $\C^\ast$-action on $\Omega^{0,1}(\End(E)) \oplus \Omega^{1,0}(\End(E))$ defined by
    \begin{eqnarray} \label{E 2}
        \xi \cdot \left( \sum_{j = 1 - \ell}^{\ell - 1} \beta_j, \sum_{j = 1 - \ell}^{\ell - 1} \varphi_j \right) = \left( \sum_{j = 1 - \ell}^{\ell - 1} \xi^j \beta_j, \sum_{j = 1 - \ell}^{\ell - 1} \xi^{j + 1} \varphi_j \right),
    \end{eqnarray}
    preserves both the negative-slice $S^{<0}_{(\bp_E, \Phi)}$ and $\Hh^1_-(\bp_E, \Phi)$.
\end{Lemma}

\begin{proof}
    First we check that the action preserves the negative-slice. Using the grading we can decompose the slice equations: the $j$th-graded piece of $D''(\beta, \varphi) + [\beta, \varphi]$ and $D'(\beta, \varphi)$ are given by
    \begin{eqnarray*}
        &&(D''(\beta, \varphi) + [\beta, \varphi])_j = \bp_E \varphi_j + [\beta_{j + 1}, \Phi] + \sum_{a + b = j} [\beta_a, \varphi_b],\\
        &&D'(\beta, \varphi)_j = \p_E \beta_j + [\Phi^\ast, \varphi_{j - 1}].
    \end{eqnarray*}
    Observe that the $\C^\ast$-action (\ref{E 2}) scales the $j$th-graded piece of the first equation by $\xi^{j + 1}$ and the $j$th-graded piece of the second equation by $\xi^j$.

    Since each $\Hh_j^1(\bp_E, \Phi) = \ker D'' \cap \ker D'$ is preserved by the action by the discussion above, $\Hh_-^1(\bp_E, \Phi)$ is also preserved by the action.
    
\end{proof}

\begin{Theorem} \label{T 1}
    If $(\bp_E, \Phi)$ is stable, the Kuranishi map gives a biholomorphism
    \begin{eqnarray*}
        k: S^{<0}_{(\bp_E, \Phi)} \xlongrightarrow{\sim} \Hh_-^1(\bp_E, \Phi).
    \end{eqnarray*}
\end{Theorem}

\begin{proof}
    The map is a homeomorphism in a neighborhood of the origin \cite[Theorem 7.3.23]{MR909698} and is clearly holomorphic. So, it is a biholomorphism in a neighborhood of the origin. Moreover, it is equivariant for the $\C^\ast$-action defined in (\ref{E 2}) since it preserves both the negative slice and $\Hh_-^1(\bp_E, \Phi)$. Hence, the local biholomorphism extends globally.
\end{proof}

The following two technical results follow closely the ideas from Proposition 3.11 in \cite{MR3949609}. The first one is the analogous version of such a result for the negative slice, and the second one produces a filtration of the $\rB \rB$-slice.

\begin{Proposition} \label{P 2}
    Suppose $(\bp_E, \Phi)$ is a stable complex variation of Hodge structure. For each
    \begin{eqnarray*}
        (\beta, \varphi) \in \Omega^{0,1}(N_-) \oplus \Omega^{1,0}(N_{<-1})
    \end{eqnarray*}
    satisfying $D''(\beta, \varphi) + [\beta, \varphi] = 0$, there exists a unique $f \in \Omega^0(N_-)$ such that the complex gauge transformation $g = \mathbf{I} + f$ takes $(\bp_E + \beta, \Phi + \varphi)$ into the slice $S^{<0}_{(\bp_E, \Phi)}$.
\end{Proposition}

\begin{proof}
    In terms of its graded pieces $f = f_{-1} + f_{-2} + \dots + f_{1- \ell}$, where $f_j \in \End_j(E)$ for $j = -1, \dots, 1 - \ell$. Observe that the complex gauge transformation $g = \mathbf{I} + f$ acts on $(\bp_E + \beta, \Phi + \varphi)$ as follows:
    \begin{equation} \label{E 5}
        \begin{array}{l}
            \bp_E + \beta \longmapsto \bp_E + \beta + \bp_E(f) +  [\beta, f] =: \bp_E + \tilde{\beta},\\
            \\
            \Phi + \varphi \longmapsto \Phi + \varphi + [\Phi, f] + [\varphi, f] =: \Phi + \tilde{\varphi},
        \end{array}
    \end{equation}
    where
    \begin{eqnarray*}
        &&\tilde{\beta} := \dsp \sum_{j=-1}^{1 - \ell} (\beta_j + \bp_E(f_j)) + \dsp \sum_{k,j= - 1}^{1 - \ell} [\beta_k, f_j] \in \Omega^{0,1}(N_-), \\
        &&\tilde{\varphi} := \dsp \sum_{j = -2}^{1-\ell} \varphi_j + \dsp \sum_{j=-1}^{1 - \ell} [\Phi, f_j] + \dsp \sum_{\substack{k=-2\\j=-1}}^{1 - \ell} [\varphi_k, f_j] \in \Omega^{1,0}(N_{<-1}),
    \end{eqnarray*}
    since the grading is additive under multiplication. To show that $(\bp_E + \beta, \Phi + \varphi) \cdot g = (\bp_E + \tilde{\beta}, \Phi + \tilde{\varphi})$ is in the negative slice for $(\bp_E, \Phi)$, we need to find $f \in \Omega^0(N_-)$ so that $(\tilde{\beta}, \tilde{\varphi})$ satisfies both slice equations:
    \begin{eqnarray*}
        &&0 = D''(\tilde{\beta}, \tilde{\varphi}) + [\tilde{\beta}, \tilde{\varphi}] = \bp_E\tilde{\varphi} + [\tilde{\beta}, \Phi] + [\tilde{\beta}, \tilde{\varphi}],\\
        &&0 = D'(\Tilde{\beta}, \tilde{\varphi}) = \p_E \tilde{\beta} + [\Phi^\ast, \tilde{\varphi}].
    \end{eqnarray*}
    Clearly $D''(\tilde{\beta}, \tilde{\varphi}) + [\tilde{\beta}, \tilde{\varphi}] = (\bp_E + \tilde{\beta})(\Phi + \tilde{\varphi}) = 0$ because $(\bp_E + \tilde{\beta}, \Phi + \tilde{\varphi})$ is a Higgs bundle. So our goal is to find $f \in \Omega^0(N_-)$ such that $D'(\tilde{\beta}, \tilde{\varphi}) = 0$. We will solve recursively for the graded pieces $f_j$ of $f$ using the $j$th-graded piece of $D'(\tilde{\beta}, \tilde{\varphi}) = 0$:
    \begin{eqnarray*}
        (\p_E \tilde{\beta} + [\Phi^\ast, \tilde{\varphi}])_j = \p_E \tilde{\beta}_j + [\Phi^\ast, \tilde{\varphi}]_j = 0.
    \end{eqnarray*}
    
    By the phrase ``mod j'' we mean any expression involving only sections $f_k$ for $k > j$ (notice that for us $j$ is negative). For example, if we write $g^{-1} = h_{-1} + h_{-2} + \dots + h_{1-\ell}$ in terms of graded pieces, then $h_j = - f_j$ mod $j$. In particular,
    \begin{eqnarray*}
        \tilde{\beta}_j = \beta_j + \bp_E(f_j) + \sum_{a + b = j} [\beta_a, f_b] = \beta_j + \bp_E(f_j) \hspace{3pt} \mbox{ mod } j,
    \end{eqnarray*}
    since $b = j - a > j$, so the sum only depends on $f_b$ for $b > j$. Moreover, if $R(f)$ is an algebraic expression in $f$ and its derivatives $\bp_Ef$, and such that the $j$th-graded piece $R(f)_j$ depends only on $f_k, k \geq j$. We will say that $R(f)_j = S(f_j)$ mod $j$, where $S$ is another such algebraic expression, if $R(f)_j - S(f_j)$ depends only on $f_k, k \geq j + 1$. For example,
    \begin{eqnarray*}
        [\Phi^\ast, \tilde{\varphi}]_j = [\Phi^\ast, \varphi]_j + [\Phi^\ast, [\Phi, f]]_j + \sum_{a + b = j - 1} [\Phi^\ast, [\varphi_a, f_b]] = [\Phi^\ast, \varphi]_j + [\Phi^\ast, [\Phi, f]]_j \hspace{3pt} \mbox{ mod } j,
    \end{eqnarray*}
    because the sum only depends on $f_b$ for $b = j - 1 - a \geq j + 1$.\\
    
    Therefore the $j$th-graded piece of $D'(\tilde{\beta}, \tilde{\varphi}) = 0$ is 
    \begin{eqnarray*}
        \p_E \beta_j + \p_E \bp_E(f_j) + [\Phi^\ast, \varphi]_j + [\Phi^\ast, [\Phi, f]]_j = 0 \hspace{3pt} \mbox{ mod } j.
    \end{eqnarray*}
    Thus, we would like to find $f$ with components $f_j$ solving the equation:
    \begin{eqnarray*}
        \p_E \bp_E(f_j) + [\Phi^\ast, [\Phi, f]]_j = - (\p_E \beta + [\Phi^\ast, \varphi])_j \hspace{3pt} \mbox{ mod } j.
    \end{eqnarray*}
    The left hand side of the equation above is $D'D''(f_j)$. Contracting with $-\sqrt-1 \Lambda$, and using the Kähler identity $(D'')^\ast = -i[\Lambda, D']$, we get
    \begin{eqnarray} \label{E 1}
        (D'')^\ast D''(f_j) = -2i\Lambda(\p_E \beta + [\Phi^\ast, \varphi])_j \hspace{3pt} \mbox{ mod } j.
    \end{eqnarray}
    Since $(\bp_E, \Phi)$ is stable, $\ker D'' = \{0\}$ on $\Omega^0(\fsl(E))$, so the operator $(D'')^\ast D''$ is invertible. Assuming that we have solved for $f_k, k > j$, the last equation gives a unique solution for $f_j$.
\end{proof}

\begin{Proposition} \label{P 3}
    Suppose $(\bp_E, \Phi)$ is a stable complex variation of Hodge structure. Let $1 \leq j_0 \leq \ell-1$ be fixed. For each
    \begin{eqnarray*}
        (\beta, \varphi) \in \Omega^{0,1}(N_{>j_0}) \oplus \Omega^{1,0}(N_{\geq j_0})
    \end{eqnarray*}
    satisfying $D''(\beta, \varphi) + [\beta, \varphi] = 0$, there is a unique $f \in \Omega^0(N_{>j_0})$ such that the complex gauge transformation $g = \mathbf{I} + f$ takes $(\bp_E + \beta, \Phi + \varphi)$ into
    \begin{eqnarray*}
        S^{\geq j_0}_{(\bp_E, \Phi)} := S^+_{(\bp_E, \Phi)} \cap ( \Omega^{0,1}(N_{>j_0}) \oplus \Omega^{1,0}(N_{\geq j_0})) \subset S^+_{(\bp_E, \Phi)}.
    \end{eqnarray*}
\end{Proposition}

\begin{proof}
    Consider $(\beta, \varphi) \in \Omega^{0,1}(N_+) \, \oplus \, \Omega^{1,0}(L \oplus N_+)$ with $(\beta_j, \varphi_{j-1}) = (0,0)$ for all $1 \leq j \leq j_0$, where $(\beta_j, \varphi_{j-1}) \in \Omega^{0,1}(\End_j(E)) \, \oplus \, \Omega^{1,0}(\End_{j - 1}(E))$. By Proposition 3.11. in \cite{MR3949609}, there is a unique $f \in \Omega^0(N_+)$ so that $g = \mathbf{I} + f$ takes $(\bp_E + \beta, \Phi + \varphi)$ into the slice $S^+_{(\bp_E, \Phi)}$. In this specific case, we can say more going back to the proof of that result in \cite{MR3949609}. When solving recursively for the graded pieces $f_j$ of $f$, observe that for $1 \leq j \leq j_0$ we need to solve the equation
    \begin{eqnarray*}
        (D'')^\ast D''(f_j) = -2i \Lambda(\p_E \beta + [\Phi^\ast, \varphi])_j \hspace{3pt} \mbox{ mod } j,
    \end{eqnarray*}
    where here ``mod $j$" means that the expression $(D'')^\ast D''(f_j) + 2i \Lambda(\p_E \beta + [\Phi^\ast, \varphi])_j$ depends only on $f_k, k \leq j - 1$. But
    \begin{eqnarray*}
        (\p_E \beta + [\Phi^\ast, \varphi])_j = \p_E \beta_j + [\Phi^\ast, \varphi_{j - 1}] = 0,
    \end{eqnarray*}
    since $(\beta_j, \varphi_{j-1}) = (0,0)$ for all $1 \leq j \leq j_0$. Thus, for $1 \leq j \leq j_0$, $f_j$ is a solution of
    \begin{eqnarray*}
        (D'')^\ast D''(f_j) = 0 \hspace{3pt} \mbox{ mod } j.
    \end{eqnarray*}
    The operator $(D'')^\ast D''$ is invertible because $\ker D'' = \{0\}$ on $\Omega^0(\fsl(E))$ since $(\bp_E, \Phi)$ is stable. Hence, $f_j = 0$ for $1 \leq j \leq j_0$, which implies that $f \in \Omega^0(N_{>j_0})$. Moreover, the complex gauge transformation $g = \mathbf{I} + f$ acts on $(\bp_E + \beta, \Phi + \varphi)$ as in (\ref{E 5}), where $(\tilde{\beta}, \tilde{\varphi}) \in \Omega^{0,1}(N_{>j_0}) \, \oplus \, \Omega^{1,0}(N_{\geq j_0})$. So we conclude that $(\tilde{\beta}, \tilde{\varphi}) \in S^{\geq j_0}_{(\bp_E, \Phi)}$.
\end{proof}

\begin{Remark} \label{R 2}
    An analogous result to Proposition \ref{P 3} holds for the negative slice. Suppose $(\bp_E + \beta, \Phi + \varphi)$ is a Higgs bundle with $(\beta, \varphi) \in \Omega^{0,1}(N_{\leq j_0}) \oplus \Omega^{1,0}(N_{< j_0})$ for some $1 - \ell \leq j_0 \leq -2$. Then equation (\ref{E 1}) becomes $(D'')^\ast D''(f_j) = 0 \, \mbox{ mod } j$ for $j_0 + 1 \leq j \leq -1 $, which implies that $f \in \Omega^0(N_{\leq j_0})$ and $g = \mathbf{I} + f$ takes $(\bp_E + \beta, \Phi + \varphi)$ into
    \begin{eqnarray*}
        S^{\leq j_0}_{(\bp_E, \Phi)} := S^{<0}_{(\bp_E, \Phi)} \cap (\Omega^{0,1}(N_{\leq j_0}) \oplus \Omega^{1,0}(N_{< j_0})) \subset S^{<0}_{(\bp_E, \Phi)}.
    \end{eqnarray*}
    Observe that with this notation $S^{<0}_{(\bp_E, \Phi)} = S^{\leq -1}_{(\bp_E, \Phi)}$.
\end{Remark}


\section{The unstable manifold at a $\VHS$ in $\Mm_\sH$} \label{DF}

In \cite{MR3949609} the authors characterized the stable manifold at a stable $\C^\ast$-fixed point, $p \in \VHS_\alpha$,
\begin{eqnarray*}
    W^{0}_\alpha(p) := \{ y \in \Mm_\sH :  \lim_{\xi \to 0} \xi \cdot y = p \},
\end{eqnarray*}
i.e. the set of points that flow down to $p$, using the $\rB \rB$-slice at that point. In particular, they proved that the restriction of (\ref{proj1}) to $S^+_{p}$ is a biholomorphism. Notice that if the stable $\C^\ast$-fixed point is a Fuchsian point from (\ref{Fuchsian}), then the map above recovers the parametrization of Hitchin sections by the affine space $\dsp \bigoplus_{j = 1}^{n - 1} H^0(K^{j+1})$. So, the result generalizes the parametrization to the fibers of every smooth stata. Moreover, there is a filtration of $W_\alpha^0$ given by the decomposition of the tangent space at $p$ into weight subspaces. Namely, combining the map (\ref{proj1}) with the filtration of the $\rB \rB$-slide from Proposition \ref{P 3}, we obtain the following:

\begin{Lemma} \label{L 3}
    There is a filtration of the Morse strata
    \begin{eqnarray*}
       W_\alpha^{0,\ell-1} \subset \dots \subset W_\alpha^{0,1} \subset W_\alpha^{0},
    \end{eqnarray*}
    where $W_\alpha^{0,j} = \dsp \bigcup_{p \in \VHS_\alpha} W_\alpha^{0, j}(p)$ and $W_\alpha^{0,j}(p) := p_\rH (S_p^{\geq j})$ for all $1 \leq j \leq \ell - 1$, for any stable $\C^\ast$-fixed point $p \in \VHS_\alpha$.
\end{Lemma}

In this section we prove similar results for the points that flow up to a stable fixed point. For a fixed point $p \in \VHS_\alpha \subset \Mm_\sH$, consider the set of Higgs bundles that flow up to it:
\begin{eqnarray*}
    W^{\infty}_\alpha(p) := \{ y \in \Mm_\sH :  \lim_{\xi \to \infty}\xi \cdot y = p \}.
\end{eqnarray*}
Notice that the limit at infinity exists if the Higgs bundles $y$ are in the nilpotent cone, i.e., $W^{\infty}_\alpha(p) \subset \NC_\alpha$. When $p$ is stable, $W^{\infty}_\alpha(p)$ is a smooth submanifold since stability is an open condition. It is called the \textit{unstable submanifold} at $p$.

\begin{Lemma} \label{L 1}
    If $[(\bp_\infty, \Phi_\infty)] \in \VHS_\alpha$ is a stable $\C^\ast$-fixed point, then $W^{\infty}_\alpha(\bp_\infty, \Phi_\infty)$ is a smooth submanifold of dimension $\dim(S^{< 0}_{(\bp_\infty, \Phi_\infty)})$.
\end{Lemma}

\begin{proof}
    By Theorem \ref{T 1}, $S^{< 0}_{(\bp_\infty, \Phi_\infty)}$ and $\Hh_-^1(\bp_\infty, \Phi_\infty)$ are biholomorphic, so they have the same dimension. Using Morse theory one can relate the grading on the tangent space at a fixed point and its unstable manifold. In \cite{MR887284}, Hitchin proved that the following function, where the adjoint is defined with respect to the harmonic metric of $(\bp_E, \Phi)$ from the nonabelian Hodge correspondence,
    \begin{eqnarray*}
        \begin{array}{rccl}
            f:&\Mm_\sH & \longrightarrow & \R^{\geq 0} \\
            &\mbox{[}(\bp_E, \Phi)\mbox{]} & \longmapsto & \|\Phi\|^2 = i \dsp \int_X \tr (\Phi \wedge \Phi^\ast) ,
        \end{array}
    \end{eqnarray*}
    is proper and, in the smooth locus, it is the moment map for the restriction of the $\C^\ast$-action to $\rU(1)$. This implies that $f$ is a Morse-Bott function on the smooth locus. Hitchin also showed \cite[Section 8]{MR1174252} that the Hessian of $f$ has eigenvalue $j$ on the subspace $\Hh_j^1(\bp_\infty, \Phi_\infty)$, since we are using the opposite grading as Hitchin, for us the Hessian of $f$ has eigenvalue $j$ on the subspace $\Hh_{-j}^1(\bp_\infty, \Phi_\infty)$. This means that
    \begin{eqnarray*}
        T_{[(\bp_\infty, \Phi_\infty)]} W^{\infty}_\alpha(\bp_\infty, \Phi_\infty) \cong \Hh_-^1(\bp_\infty, \Phi_\infty).
    \end{eqnarray*}
    Thus, $\dim(W^{\infty}_\alpha(\bp_\infty, \Phi_\infty)) = \dim(\Hh_-^1(\bp_\infty, \Phi_\infty)) = \dim(S^{< 0}_{(\bp_\infty, \Phi_\infty)})$.
\end{proof}

\begin{Proposition} \label{P 1}
    Let $[(\bp_\infty, \Phi_\infty)] \in \VHS_\alpha$ be stable, and $[(\bp_E, \Phi)] \in \NC$. Then $[(\bp_E, \Phi)] \in W^{\infty}_\alpha(\bp_\infty, \Phi_\infty)$ if and only if after a gauge transformation
    \begin{eqnarray*}
        \Phi - \Phi_\infty \in \Omega^{1,0}(N_{<-1}) \hspace{10pt} \mbox{ and } \hspace{10pt} \bp_E - \bp_\infty \in \Omega^{0,1}(N_-). 
    \end{eqnarray*}
\end{Proposition}

\begin{proof}
    Suppose $\varphi = \Phi - \Phi_\infty \in \Omega^{1,0}(N_{<-1})$ and $\beta = \bp_E - \bp_\infty \in \Omega^{0,1}(N_-)$. Using the grading defined on $\End(E)$ we can write
    \begin{eqnarray*}
        \varphi = \sum_{j= -2}^{1 - \ell} \varphi_j \hspace{10pt} \mbox{ and } \hspace{10pt} \beta = \sum_{j = -1}^{1 - \ell} \beta_j,
    \end{eqnarray*}
    where $\varphi_j \in \Omega^{1,0}(\End_j(E))$ and $\beta_j \in \Omega^{0,1}(\End_j(E))$ for any $1 - \ell \leq j \leq -1$. Let $\xi \in \C^\ast$ and consider the gauge transformation
    \begin{eqnarray} \label{gauge}
        g_\xi = \begin{pmatrix}
            \xi^a \mathbf{I}_{E_1} & &  \\
            & \xi^{a + 1} \mathbf{I}_{E_2} & &  \\
            & & \ddots &  \\
            & & & \xi^{a + \ell - 1} \mathbf{I}_{E_\ell}
        \end{pmatrix},
    \end{eqnarray}
    with $\displaystyle \sum_{j = 0}^{\ell - 1} \rk(E_{j+1})(a + j) = 0$ because we want $\det(g_\xi) = 1$. Then
    \begin{eqnarray*}
        (\bp_E, \xi \Phi) \cdot g_\xi = (\bp_\infty + \beta, \xi \Phi_\infty + \xi \varphi) \cdot g_\xi = \left( \bp_\infty + \sum_{j = -1}^{1 - \ell} \xi^j \beta_j, \Phi_\infty + \sum_{j = -2}^{1 - \ell} \xi^{j+1} \varphi_j \right),
    \end{eqnarray*}
    and $\displaystyle \lim_{\xi \to \infty}[(\bp_E, \xi \Phi)] = [(\bp_\infty, \Phi_\infty)]$.\\
    
    For the other direction, let $[(\bp_E, \Phi)] \in W^{\infty}_\alpha(\bp_\infty, \Phi_\infty)$. Since $(\bp_\infty, \Phi_\infty)$ is stable, the map $p_\rH$ from (\ref{proj1}) is a diffeomorphism from a neighborhood of the origin in $S_{(\bp_\infty, \Phi_\infty)}$ to an open neighborhood $\mathcal{U}  \subset \Mm_\sH$ of $[(\bp_\infty, \Phi_\infty)]$. Then, there exists $\xi \in \C^\ast$ such that $[(\bp_E, \xi \Phi)] \in \mathcal{U}$ because $\displaystyle \lim_{\xi \to \infty}[(\bp_E, \xi \Phi)] = [(\bp_\infty, \Phi_\infty)]$. Observe that $[(\bp_E, \xi \Phi)] \in \NC$ since so is $[(\bp_E, \Phi)]$, and $\displaystyle \lim_{\eta \to \infty}[(\bp_E, \eta (\xi \Phi))] = [(\bp_\infty, \Phi_\infty)]$. This means that $[(\bp_E, \xi \Phi)] \in W^{\infty}_\alpha(\bp_\infty, \Phi_\infty)$. Therefore,
    \begin{eqnarray*}
        [(\bp_E, \xi \Phi)] \in W^{\infty}_\alpha(\bp_\infty, \Phi_\infty) \cap \mathcal{U}.
    \end{eqnarray*}
    
    Using the smooth splitting $E = E_1 \oplus \dots \oplus E_\ell$ we can decompose $\beta \in \Omega^{0,1}(\End(E))$ and $\varphi \in \Omega^{1,0}(\End(E))$ as $\beta = \dsp \sum_{j = \ell - 1}^{1 - \ell} \beta_j$ and $\varphi = \dsp \sum_{j = \ell - 1}^{1 - \ell} \varphi_j$. Then, we can write the negative slice as
    \begin{eqnarray*}
        S_{(\bp_\infty, \Phi_\infty)}^{< 0} = \left\{ (\beta, \varphi) \in S_{(\bp_\infty, \Phi_\infty)} : \beta = \dsp \sum_{j = -1}^{1 - \ell} \beta_j \; \mbox{ and } \; \varphi = \dsp \sum_{j = -2}^{1 - \ell} \varphi_j \right\}. 
    \end{eqnarray*}
    
    Notice that
    \begin{eqnarray*}
        W_\alpha^\infty(\bp_\infty, \Phi_\infty) \cap p_{\rH}(S_{(\bp_\infty, \Phi_\infty)}) \subset p_{\rH}(S_{(\bp_\infty, \Phi_\infty)}^{< 0}),
    \end{eqnarray*}
    by the argument in the previous implication. We know that $\dim(W^{\infty}_\alpha(\bp_\infty, \Phi_\infty)) = \dim(S^{< 0}_{(\bp_\infty, \Phi_\infty)})$ by Lemma \ref{L 1}, and the above inclusion is an equality.

    For each $[(\bp_E, \Phi)] \in W^{\infty}_\alpha(\bp_\infty, \Phi_\infty)$ there is $\xi \in \C^\ast$ so that $[(\bp_E, \xi \Phi)] \in \mathcal{U} = p_\rH(S_{(\bp_\infty, \Phi_\infty)})$, this means that $[(\bp_E, \xi \Phi)] \in p_\rH(S_{(\bp_\infty, \Phi_\infty)}^{< 0})$. Therefore there is a representative $(\bp_E, \xi \Phi)$ of $[(\bp_E, \xi \Phi)]$ such that
    \begin{eqnarray*}
        (\bp_E, \xi \Phi) = \left( \bp_\infty + \sum_{j = -1}^{1 - \ell} \beta_j, \Phi_\infty + \sum_{j = -2}^{1 - \ell} \varphi_j \right).
    \end{eqnarray*}
    A representative of $[(\bp_E, \Phi)] = \xi^{-1} \cdot [(\bp_E, \xi \Phi)]$ is given by
    \begin{eqnarray*}
        (\bp_E, \Phi) = \left( \bp_\infty + \sum_{j = -1}^{1 - \ell} \beta_j, \xi^{-1} \Phi_\infty + \sum_{j = -2}^{1 - \ell} \xi^{-1} \varphi_j \right).
    \end{eqnarray*}
    Using the inverse of the holomorphic gauge transformation $g_\xi$ from (\ref{gauge}) we get
    \begin{eqnarray*}
        (\bp_E, \xi \Phi) \cdot g_\xi^{-1} &=& \left( \bp_\infty + \sum_{j = -1}^{1 - \ell} \beta_j, \xi^{-1} \Phi_\infty + \sum_{j = -2}^{1 - \ell} \xi^{-1} \varphi_j \right) \cdot g_\xi^{-1} \\
        &=& \left( \bp_\infty + \sum_{j = -1}^{1 - \ell} \xi^{-j} \beta_j, \Phi_\infty + \sum_{j = -2}^{1 - \ell} \xi^{-j -1} \varphi_j \right).
    \end{eqnarray*}
    Thus, after a gauge transformation:
    \begin{eqnarray*}
        \bp_E - \bp_\infty \in \Omega^{0,1}(N_-) \hspace{10pt} \mbox{ and } \hspace{10pt} \Phi - \Phi_\infty \in \Omega^{1,0}(N_{<-1}).
    \end{eqnarray*}
    Hence, every $[(\bp_E, \Phi)] \in W^{\infty}_\alpha(\bp_\infty, \Phi_\infty)$ is gauge equivalent to $(\bp_E, \Phi)$ with $\bp_E - \bp_\infty \in \Omega^{0,1}(N_-)$ and $\Phi - \Phi_\infty \in \Omega^{1,0}(N_{<-1})$.
\end{proof}

\begin{Corollary} \label{C 1}
    Let $[(\bp_E, \Phi)] \in \VHS_\alpha$ be stable. Then the map
    \begin{eqnarray*}
        \begin{array}{rccl}
            p_\rH : & S^{<0}_{(\bp_E, \Phi)} & \lra & W^{\infty}_\alpha(\bp_E, \Phi)   \\
            &(\beta, \varphi) & \longmapsto & [(\bp_E + \beta, \Phi + \varphi)],
        \end{array}
    \end{eqnarray*}
    is a biholomorphism.
\end{Corollary}

\begin{proof}
    First we show that the map is injective. Let $(\beta_1, \varphi_1), (\beta_2, \varphi_2) \in S^{<0}_{(\bp_E, \Phi)}$ be such that for some $g \in \Gg(E)$,
    \begin{eqnarray*}
        (\bp_E + \beta_1, \Phi + \varphi_1) \cdot g = (\bp_E + \beta_2, \Phi + \varphi_2).
    \end{eqnarray*}
    Let $\xi \in \C^\ast$, scaling the previous equality we obtain
    \begin{eqnarray*}
        (\bp_E + \beta_1, \xi \Phi + \xi \varphi_1) \cdot g = (\bp_E + \beta_2, \xi \Phi + \xi \varphi_2).
    \end{eqnarray*}
    Now applying the gauge transformation $g_\xi$ from (\ref{gauge}), we have
    \begin{eqnarray*}
        \left( \bp_E + \sum_{j = -1}^{1 - \ell} \xi^{j} \beta_{1,j}, \Phi + \sum_{j = -2}^{1 - \ell} \xi^{j + 1} \varphi_{1,j} \right) \cdot (g_\xi^{-1} \circ g \circ g_\xi) = \left( \bp_E + \sum_{j = -1}^{1 - \ell} \xi^{j} \beta_{2,j}, \Phi + \sum_{j = -2}^{1 - \ell} \xi^{j + 1} \varphi_{2,j} \right).
    \end{eqnarray*}
    For $\xi$ sufficiently large, this is a gauge equivalence between points in a neighborhood of the origin in $S^{<0}_{(\bp_E, \Phi)}$. The map $p_\rH$ is injective there by Theorem 7.3.17 in \cite{MR909698}, so $(\beta_1, \varphi_1) = (\beta_2, \varphi_2)$, which proves that $p_\rH$ is injective on all of $S^{<0}_{(\bp_E, \Phi)}$. To show that $p_\rH$ is surjective let $[(\bp_E', \Phi')] \in W^{\infty}_\alpha(\bp_E, \Phi)$. Then, by Proposition \ref{P 1}, after a gauge transformation we have
    \begin{eqnarray*}
        \varphi := \Phi' - \Phi \in \Omega^{(1,0)}(N_{<-1}) \hspace{10pt} \mbox{ and } \hspace{10pt} \beta := \bp_E' - \bp_E \in \Omega^{(0,1)}(N_-).
    \end{eqnarray*}
    So we can write the original Higgs bundle as $[(\bp_E', \Phi')] = [(\bp_E + \beta, \Phi + \varphi)]$. We need to check if $(\beta, \varphi)$ is in $S^{<0}_{(\bp_E, \Phi)}$. Notice that
    \begin{eqnarray*}
        D''(\beta, \varphi) + [\beta, \varphi] 
        = \bp_E \varphi + [\beta, \Phi'] = \bp_E' \Phi' = 0.
    \end{eqnarray*}
    Applying Proposition \ref{P 2}, we know that after a gauge transformation $(\beta, \varphi) \in S^{<0}_{(\bp_E, \Phi)}$. Thus, $p_\rH$ is surjective. Moreover, $p_\rH$ is holomophic, which concludes the proof.
\end{proof}

As a consequence of this last result, the filtration we obtained for the negative slice at a stable fixed point in Remark \ref{R 2} produces a filtration of $W_\alpha^\infty$. 

\begin{Corollary} \label{C 2}
    There is a filtration of $W_\alpha^\infty$ given by
    \begin{eqnarray*}
       W_\alpha^{\infty,1-\ell} \subset \dots \subset W_\alpha^{\infty,-2} \subset W_\alpha^{\infty},
    \end{eqnarray*}
    where $W_\alpha^{\infty,j} = \dsp \bigcup_{p \in \VHS_\alpha} W_\alpha^{\infty,j}(p)$ and $W_\alpha^{\infty,j}(p) := p_\rH (S_p^{\leq j})$ for all $1 - \ell \leq j \leq - 2$, and $W_\alpha^{\infty}(p) = p_\sH(S_p^{< 0})$, for each stable $\C^\ast$-fixed point $p \in \VHS_\alpha$.
\end{Corollary}

\begin{Remark}
    Recall that $\NC_\alpha = W_\alpha^\infty$, which means that Corollary \ref{C 2} gives a filtration of the irreducible components of the nilpotent cone:
    \begin{eqnarray*}
        \NC_\alpha^{1-\ell} \subset \dots \subset \NC_\alpha^{-2} \subset \NC_\alpha^{-1} = \NC_\alpha.
    \end{eqnarray*}
\end{Remark}

\section{Upward conformal limits} \label{S 5}

In this section we use the results from Sections \ref{S 3} and \ref{DF} to prove Theorems \ref{MT} and \ref{MT2}, and Corollary \ref{MT3}.

\subsection{Existence of the upward conformal limit}

Let $[(\bp_E, \Phi)] \in \VHS_\alpha$ be a stable Higgs bundle, and let $h$ be the harmonic metric for $(\bp_E, \Phi)$. By Corollary \ref{C 1}, any point in $W^\infty_\alpha(\bp_E, \Phi)$ is of the form $[(\bp_E + \beta, \Phi + \varphi)]$ for some
$\textbf{u}:= (\beta, \varphi) \in S_{(\p_E, \Phi)}^{< 0}$. Set
\begin{eqnarray*}
    \bp_\textbf{u} := \bp_E + \beta \hspace{10pt} \mbox{ and } \hspace{10pt} \Phi_\textbf{u} := \Phi + \varphi.
\end{eqnarray*}
Given $R > 0$, we can define a family of flat connections parametrized by $R$ as follows. Scaling the Higgs field, we obtain a collection of stable Higgs bundles $(\bp_\textbf{u}, R \Phi_\textbf{u})$ and their harmonic metrics $h(\textbf{u},R)$, that make the connections $\bp_\textbf{u} + \p_\textbf{u}^{h(\textbf{u}, R)} + R\Phi_\textbf{u} + R \Phi_\textbf{u}^{\ast_{h(\textbf{u}, R)}}$ flat. For $\xi \in \C^\ast$, the connections
\begin{eqnarray*}
    \bp_\textbf{u} + \p_\textbf{u}^{h(\textbf{u}, R)} + \xi^{-1} R\Phi_\textbf{u} + \xi R \Phi_\textbf{u}^{\ast_{h(\textbf{u}, R)}} 
\end{eqnarray*}
are also flat. If we keep $\hbar := \xi R^{-1}$ fixed, we obtain a family of flat connections parametrized by $R$:
\begin{eqnarray} \label{E FC}
    D_{(\textbf{u}, R)} := \bp_\textbf{u} + \p_\textbf{u}^{h(\textbf{u}, R)} + \hbar^{-1} \Phi_\textbf{u} + \hbar R^2 \Phi_\textbf{u}^{\ast_{h(\textbf{u}, R)}}. 
\end{eqnarray}

\begin{Theorem} \label{CL}
    Let $[(\bp_E, \Phi)] \in \VHS_\alpha$ be stable and $\textbf{u}:= (\beta, \varphi) \in S_{(\p_E, \Phi)}^{< 0}$. As $R \to \infty$, the family of flat connections from (\ref{E FC}) converges smoothly to
    \begin{eqnarray*}
        \bp_\textbf{u} + \p_E + \hbar^{-1} \Phi_\textbf{u} + \hbar \Phi^\ast.
    \end{eqnarray*}
    We call this limiting flat connection the $\hbar$-upward conformal limit of $[(\bp_\textbf{u}, \Phi_\textbf{u})]$, and we denote it by $CL_{\hbar, \infty}([(\bp_\textbf{u}, \Phi_\textbf{u})])$. In particular,
    \begin{eqnarray*}
        CL_{\hbar, \infty}([(\bp_\textbf{u}, \Phi_\textbf{u})]) = \hbar^{-1} \cdot p_{\hbar}(\textbf{u}).
    \end{eqnarray*}
\end{Theorem}

Restating the last result we obtain Theorem \ref{MT}.
\begin{Corollary}
    Let $[(\bp_\infty, \Phi_\infty)] \in \VHS_\alpha$ be stable and $[(\bp_E, \Phi)] \in W^\infty_\alpha(\bp_\infty, \Phi_\infty)$. Then the $\hbar$-upward conformal limit of $[(\bp_E, \Phi)]$ exists.
\end{Corollary}

In order to prove the existence of the $\hbar$-upward conformal limit we will follow the same philosophy as \cite[Proposition 5.1]{MR3949609} and \cite[Theorem 3.2]{MR4214341}. To make the proof more comprehensible, we break it down into the next four lemmas. By direct computation one has the following:

\begin{Lemma}
    Let $g \in \Gg(E)$. Define a new metric $g(h)$ by $\ip{u}{v}_{g(h)} := \ip{gu}{gv}_h$. Let $k := g^* g$. Then we have for any Higgs bundle $(\bp_E, \Phi)$,
    \begin{equation} \label{1}
        \Phi^{*_{g(h)}} = k^{-1} \Phi^{\ast_h} k,
    \end{equation}
    \begin{equation} \label{2}
        \p^{g(h)}_E = k^{-1} \circ \p_E \circ k,
    \end{equation}
    \begin{equation} \label{3}
        F_{(\bp_E, g(h))} = F_{(\bp_E, h)} + \bp_E(k^{-1} \p_E(k)),
    \end{equation}
where $F_{(\bp_E, h')}$ denotes the curvature of the Chern connection associated to $\bp_E$ and a metric $h'$.
\end{Lemma}

\begin{Lemma} \label{LCL1}
    Let $[(\bp_E, \Phi)] \in \VHS_\alpha$ be stable and $\textbf{u}:= (\beta, \varphi) \in S_{(\p_E, \Phi)}^{< 0}$. Let $R > 0$ and consider the metric $h_R' = g_R(h)$ where $h$ be the harmonic metric for $(\bp_E, \Phi)$ and
    \begin{eqnarray} \label{Gauge}
        g_R := \begin{pmatrix}
            R^{\frac{m_1}{2}} & & 0\\
            & \ddots & &\\
            0 & & R^{\frac{m_\ell}{2}}
        \end{pmatrix}, \hspace{10pt} k_R := \begin{pmatrix}
            R^{m_1} & & 0\\
            & \ddots & &\\
            0 & & R^{m_\ell}
        \end{pmatrix},
    \end{eqnarray}
    with $m_j - m_{j + 1} = 2, j = 1, \dots, \ell - 1$, and $\dsp \sum_{i = 1}^\ell (\rk E_j) m_j = 0$. Then, the family of connections
    \begin{eqnarray*}
        \tilde{D}_{(\textbf{u}, R)} := \bp_\textbf{u} + \p_\textbf{u}^{h_R'} + \hbar^{-1} \Phi_\textbf{u} + \hbar R^2 \Phi_\textbf{u}^{\ast_{h_R'}}
    \end{eqnarray*}
    can be written as
    \begin{eqnarray*}
        \bp_\textbf{u} + \p_E - \dsp\sum_{j = -1}^{1 - \ell} R^{2 j} (\beta_j)^\ast + \hbar^{-1} \Phi_\textbf{u} + \hbar \Phi^\ast + \hbar \dsp\sum_{j = -2}^{1 - \ell} R^{2j + 2} (\varphi_j)^\ast.
    \end{eqnarray*}
\end{Lemma}

\begin{proof}
    Notice that the adjoint of the Higgs field with respect to the metric $h_R'$, and the Chern connection associated to $\bp_\textbf{u}$ and the metric $h_R'$ are
    \begin{eqnarray*}
        \Phi_\textbf{u}^{\ast_{h_R'}} = \Phi^{\ast_{h_R'}} + \varphi^{\ast_{h_R'}} \hspace{10pt} \mbox{ and } \hspace{10 pt} \bp_\textbf{u} + \p_\textbf{u}^{h_R'} = \bp_\textbf{u} + \p_E^{h_R'} - \beta^{\ast_{h_R'}},
    \end{eqnarray*}
    since $\p_E - \beta^\ast$ is the Chern connection associated to $\bp_\textbf{u}$ and the metric $h$. Therefore, we need to compute $\Phi^{\ast_{h_R'}}, \p_E^{h_R'},  \varphi^{\ast_{h_R'}}$ and $\beta^{\ast_{h_R'}}$ using (\ref{1}) and (\ref{2}). Keeping in mind that $m_j - m_{j + 1} = 2, j = 1, \dots, \ell - 1$, we have
    \begin{eqnarray*}
        \Phi^{\ast_{h_R'}} = k_R^{-1} \Phi^{\ast_h} k_R = 
        R^{-2} \Phi^\ast, \hspace{10 pt} \mbox{ and } \hspace{10 pt} \p_E^{h_R'} = k_R^{-1} \circ \p_E \circ k_R = \p_E.
    \end{eqnarray*}
    For $\varphi^{\ast_{h_R'}}$ and $\beta^{\ast_{h_R'}}$, recall that the negative slice $S_{(\bp_E, \Phi)}^{< 0} \subset S_{(\bp_E, \Phi)}$ is
    \begin{eqnarray*}
        S_{(\bp_E, \Phi)}^{< 0} = \left\{ (\beta, \varphi) \in S_{(\bp_E, \Phi)} | \beta = \dsp \sum_{j = -1}^{1 - \ell} \beta_j \mbox{ and } \varphi = \dsp \sum_{j = - 2}^{1 - \ell} \varphi_j \right\},
    \end{eqnarray*}
    where $\beta_j \in \Omega^{0,1}(\End_j(E))$ and $\varphi_j \in \Omega^{1,0}(\End_j(E))$ (in particular, notice that they are lower triangular matrices). Then
    \begin{eqnarray*}
        \varphi^{\ast_{h_R'}} = \dsp \sum_{j = - 2}^{1 - \ell} (\varphi_j)^{\ast_{h_R'}} = \dsp \sum_{j= -2}^{1 - \ell} k_R^{-1} (\varphi_j)^\ast k_R = \dsp \sum_{j = -2}^{1 - \ell} R^{2j} (\varphi_j)^\ast,
    \end{eqnarray*}
    and similarly for $\beta^{\ast_{h_R'}}$. Thus,
    \begin{eqnarray*}
        \tilde{D}_{(\textbf{u}, R)} = \bp_\textbf{u} + \p_E - \dsp\sum_{j = -1}^{1 - \ell} R^{2 j} (\beta_j)^\ast + \hbar^{-1} \Phi_\textbf{u} + \hbar \Phi^\ast + \hbar \dsp\sum_{j = -2}^{1 - \ell} R^{2j + 2} (\varphi_j)^\ast.
    \end{eqnarray*}
\end{proof}

\begin{Lemma} \label{LCL2}
    Let $[(\bp_E, \Phi)] \in \VHS_\alpha$ be stable and  $\textbf{u}:= (\beta, \varphi) \in S_{(\p_E, \Phi)}^{< 0}$. Denote by $h$ and $h(\textbf{u}, R)$ the harmonic metrics for $(\bp_E, \Phi)$ and $(\bp_\textbf{u}, R\Phi_\textbf{u})$, respectively. Then, the harmonic metrics $h(\textbf{u}, R)$ can be written uniquely as
    \begin{eqnarray*}
        h(\textbf{u}, R) = g(h_R'),
    \end{eqnarray*}
    where $h_R' := g_R(h)$ as in Lemma \ref{LCL1}, and $g := \exp(f_R)$ with $f_R$ a traceless $h_R'$-hermitian endomorphism.
\end{Lemma}

\begin{proof}
    The gauge group $\Gg(E)$ acts transitively in the space of hermitian metrics on $E$, i.e., there exists $g \in \Gg(E)$ so that $h(\textbf{u}, R) = g(h_R')$. Observe that $(\bp_E, \Phi) \cdot g_R = (\bp_E, R \Phi)$ and the metric $h_R'$ is the harmonic metric for $(\bp_E, R \Phi)$. 
    Moreover, one checks that $(\bp_\textbf{u}, R \Phi_\textbf{u}) = (\bp_E + \beta, R \Phi + R \varphi)$ with $(\beta, R \varphi) \in S^{< 0}_{(\bp_E , R \Phi)}$. 
    Also, notice that in the Higgs bundle moduli space we have
    \begin{eqnarray*}
        \lim_{R \to \infty} (\bp_\textbf{u}, R \Phi_\textbf{u}) \cdot g_R^{-1} 
        = \lim_{R \to \infty} \dsp \left(\bp_E + \sum_{j = -1}^{1 - \ell} R^j \beta_j, \Phi + \sum_{j = -2}^{1 - \ell} R^{j + 1} \varphi_j \right) = (\bp_E, \Phi),
    \end{eqnarray*}
    which means that the scaling we are doing on $(\bp_\textbf{u}, \Phi_\textbf{u})$ is asymptotic to the gauge orbit of $(\bp_E, \Phi)$ for $R > 0$ large (see Figure \ref{fig1}).
    \begin{figure}[h!]
        \centering
        \includegraphics[width=0.5\textwidth]{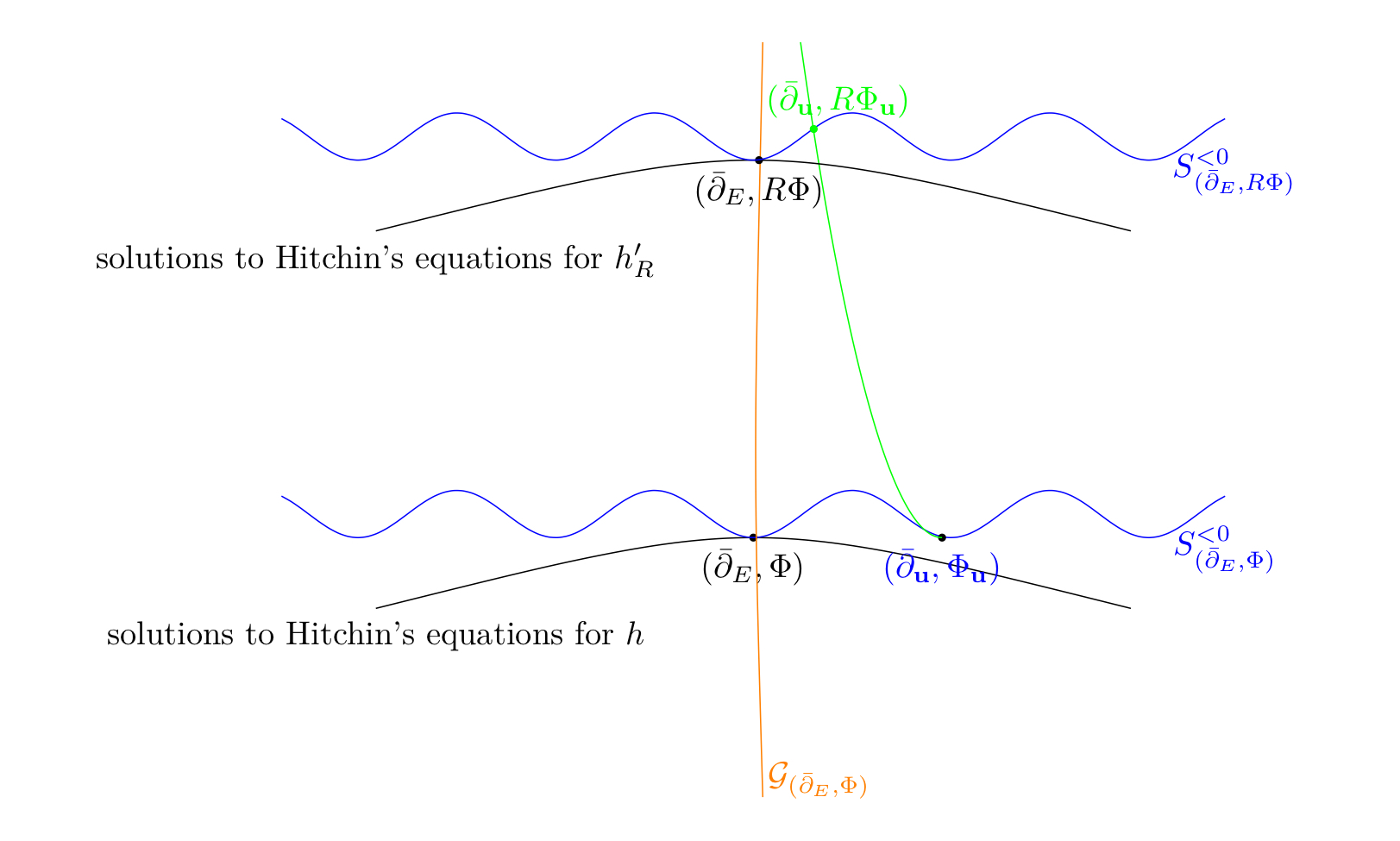}
        \caption{Asymptotic behavior of the scaling.}
        \label{fig1}
    \end{figure}
    
    Thus, by the implicit function theorem, there is a unique $f_R$ a traceless $h_R'$-hermitian endomorphism so that $g = \exp(f_R)$.
\end{proof}

The next technical result is analogous to Lemma 3.3 in \cite{MR4214341}.

\begin{Lemma} \label{LCL3}
    The $f_R$ traceless $h_R'$-hermitian endomorphism from Lemma \ref{LCL2} satisfies $f_R \to 0$ as $R \to \infty$.
\end{Lemma}

\begin{proof}
    Recall that $h(\textbf{u}, R)$ is a solution of Hitchin's equations for $(\bp_\textbf{u}, R \Phi_\textbf{u})$, i.e.,
    \begin{eqnarray*}
        0 &=& F_{(\bp_\textbf{u}, h(\textbf{u},R))} + [R \Phi_\textbf{u}, (R \Phi_\textbf{u})^{\ast_{h(\textbf{u},R)}}]\\
        &=& F_{(\bp_\textbf{u}, h_R')} + \bp_\textbf{u}(e^{-2 f_R} \p_\textbf{u}^{h_R'}(e^{2 f_R})) + [\Phi_\textbf{u}, e^{-2 f_R} (R^2 \Phi_\textbf{u}^{\ast_{h_R'}}) e^{2 f_R}]\\
        &=:& N_{(\textbf{u},R)}(f_R),
    \end{eqnarray*}
    where we used (\ref{1}), (\ref{3}), and $(R \Phi_\textbf{u})^{\ast_{h(\textbf{u},R)}} = R \Phi_\textbf{u}^{\ast_{h(\textbf{u},R)}}$. For any fixed $R$, this gives us a non-linear operator
    \begin{eqnarray*}
        N_{(\textbf{u},R)} : \dsp \Omega^0(\End_{h_R'}(E)) \lra \Omega^2(\End_{h_R'}(E)),
    \end{eqnarray*}
    where $\End_{h_R'}(E)$ is the subalgebra of traceless $h_R'$-hermitian endomorphisms.
    Note that for any $R > 0$ and $f \in \End_{h_R'}(E)$, we have $N_{(\textbf{u}, R)}(f) = 0$ if and only if $g(h_R')$ with $g = \exp(f)$, is the harmonic metric for $(\bp_\textbf{u}, R \Phi_\textbf{u})$. Since $h(\textbf{u}, R) = g(h_R')$ with $g = \exp(f_R)$, our goal is to find such an endomorphism.

    We proceed in steps.
    \begin{itemize}
        \item Step 1: $N_{(\textbf{u}, \infty)}(0) = 0$.
        Let us record 
        \begin{eqnarray} \label{5}
            F_{(\bp_\textbf{u}, h_R')} 
            = F_{(\bp_E, h)} + \dsp \left[\beta, \sum_{j= -1}^{1 - \ell} R^{2 j}(\beta_j)^\ast\right] + \p_E \beta - \sum_{j= -1}^{1 - \ell} R^{2 j} \bp_E (\beta_j)^\ast,
        \end{eqnarray}
        and
        \begin{eqnarray} \label{6}
            R^2 \Phi_\textbf{u}^{\ast_{h_R'}} = \Phi^\ast + \sum_{j = -2}^{1 - \ell} R^{2 j + 2}(\varphi_j)^\ast,
        \end{eqnarray}
        where the adjoints with respect to the metric $h_R'$ were computed in the proof of Lemma \ref{LCL1}. Then, by (\ref{5}) and (\ref{6}) we have
        \begin{eqnarray*}
            N_{(\textbf{u}, \infty)}(0) &=& \dsp \lim_{R \to \infty} N_{(\textbf{u}, R)}(0) = \lim_{R \to \infty} \lb F_{(\bp_\textbf{u}, h_R')} + [\Phi_\textbf{u}, R^2 \Phi_\textbf{u}^{\ast_{h_R'}}] \rb\\
            &=& F_{(\bp_E, h)} + \p_E \beta + [\Phi_\textbf{u}, \Phi^\ast] 
            = 0,
        \end{eqnarray*}
        since $h$ is the harmonic metric for $(\bp_E, \Phi)$ and $(\beta, \varphi) \in S^{< 0}_{(\bp_E, \Phi)}$.

        \item Step 2: compute $d N_{(\textbf{u}, \infty)}(0)$. Recall that $f_R$ is a (traceless) $h_R'$-hermitian endomorphism. So, $f_R$ is determined by a section $(f_L(R), f_-(R))$ of $L \oplus N_-$, i.e.,
        \begin{eqnarray*}
            f_R = \dsp f_L(R) + f_-(R) + \sum_{j = -1}^{1 - \ell} R^{2 j}(f_{-,j}(R))^\ast,
        \end{eqnarray*}
        where the last term comes from $f_R \in \End_{h_R'}(E)$. Consider the linearization of $N(\textbf{u}, R)$ in $R$ at $R = \infty$. Changing variables we consider $f_{\frac{1}{R}}$, and we denote its derivatives at $\frac{1}{R} = 0$ by $\dot{f}$. 
        \begin{eqnarray*}
            \frac{1}{2} d N_{(\textbf{u}, R)}(f_{\frac{1}{R}}) \dot{f} &=& \bp_\textbf{u}\left(e^{-2 f_{\frac{1}{R}}} \p_\textbf{u}^{h_R'} e^{2 f_{\frac{1}{R}}} \dot{f} - \dot{f} e^{-2 f_{\frac{1}{R}}} \p_\textbf{u}^{h_R'} e^{2 f_{\frac{1}{R}}}\right)\\
            &+& \left[\Phi_\textbf{u}, e^{-2 f_{\frac{1}{R}}} (R^2 \Phi_\textbf{u}^{\ast_{h_R'}}) e^{2 f_{\frac{1}{R}}} \dot{f}
            - \dot{f} e^{-2 f_{\frac{1}{R}}} (R^2 \Phi_\textbf{u}^{\ast_{h_R'}}) e^{2 f_{\frac{1}{R}}}\right].
        \end{eqnarray*}
        Evaluating at $\frac{1}{R} = 0$, we have
        \begin{equation} \label{7}
            \frac{1}{2} d N_{(\textbf{u}, \infty)}(0) \dot{f} = \bp_\textbf{u} \p_E \dot{f} + [\Phi_\textbf{u}, [\Phi^\ast, \dot{f}]],
        \end{equation}
        since as $R \to \infty$ we know $\p_\textbf{u}^{h_R'} \to \p_E$ and $R^2 \Phi_\textbf{u}^{\ast_{h_R'}} \to \Phi^\ast$ by (\ref{6}).\\

        \item Step 3: $d N_{(\textbf{u}, \infty)}(0)$ is injective, i.e., $\ker d N_{(\textbf{u}, \infty)}(0) = \{0\}$. Suppose $\dot{f} \in \ker d N_{(\textbf{u}, \infty)}(0)$. We can write the decomposition into components as $\dot{f} = \dot{f}_L + \dot{f}_-$ (observe that the other term vanishes because $\frac{1}{R} = 0$). Also recall Hitchin's equations
        \begin{eqnarray*}
            F_{(\bp_E, h)} + [\Phi, \Phi^\ast] = \bp_E \p_E + \p_E \bp_E + \Phi \Phi^\ast + \Phi^\ast \Phi = 0.
        \end{eqnarray*}
        Then, by (\ref{7}) and the previous expression (and the Jacobi identity), we get the following
        \begin{eqnarray} \label{8}
            0 &=& \frac{1}{2} d N_{(\textbf{u}, \infty)}(0) \dot{f} = \bp_\textbf{u} \p_E \dot{f} + [\Phi_\textbf{u}, [\Phi^\ast, \dot{f}]]\\ \nonumber 
            &=& - \p_E \bp_E \dot{f} - [\Phi^\ast, [\Phi, \dot{f}]] \color{black} + [\beta, \p_E \dot{f}] + [\varphi, [\Phi^\ast, \dot{f}]].
        \end{eqnarray}
        Notice that $[\beta, \p_E \dot{f}] + [\varphi, [\Phi^\ast \dot{f}]]$ does not preserve the grading, in particular, it does not contain any diagonal terms. So we can just look at $\dot{f}_L$, i.e.,
        \begin{eqnarray*}
            0 = \p_E \bp_E \dot{f}_L + [\Phi^\ast, [\Phi, \dot{f}_L]].
        \end{eqnarray*}
        Contracting with $i \Lambda$ and using $(D'')^\ast = -i [\Lambda, D']$ we have $(D'')^\ast D''(\dot{f}_L) = 0$. 
        By the assumption of stability, $\ker D'' = \{0\}$ on $\Omega^0(\fsl(E))$. Thus, $\dot{f}_L = 0$ and we can write $\dsp \dot{f} = \dot{f}_- = \dsp \sum_{j = -1}^{1 - \ell} \dot{f}_{-,j}$. Using equation (\ref{8}) again and that $[\beta, \p_E \dot{f}] + [\varphi, [\Phi^\ast \dot{f}]]$ does not preserve the grading, we have that $(D'')^\ast D'' (\dot{f}_{-,-1}) = 0$. By the same argument as before, $\dot{f}_{-,-1} = 0$. Repeating this process we have $\dot{f}_{-, k} = 0$ for $j < k < -1$, which implies $(D'')^\ast D''(\dot{f}_{-,j}) = 0$. Hence, $\dot{f}_- = 0$ and $\dot{f} = 0$. So, $d N_{(\textbf{u}, \infty)}(0)$ is injective.

        \item Step 4: $d N_{(\textbf{u}, \infty)}(0)$ is surjective. 
        Recall that the index of an operator is independent of deformations through elliptic operators. Deforming $d N_{(\textbf{u}, \infty)}(0)$ among elliptic operators to $\bp_E \p_E$, which is self-adjoint (i.e., $\ker(\bp_E \p_E) \cong \mbox{coker}(\bp_E \p_E)$), we have
        \begin{eqnarray*}
            \ind(d N_{(\textbf{u}, \infty)}(0)) = \ind (\bp_E \p_E) = 0.
        \end{eqnarray*}
        By the Atiyah-Singer Theorem, we have that $\dim \mbox{coker} d N_{(\textbf{u}, \infty)}(0) = 0$, since $\ker d N_{(\textbf{u}, \infty)}$ is injective. This means that $d N_{(\textbf{u}, \infty)}(0)$ is surjective.
    \end{itemize}

    We showed that $d N_{(\textbf{u}, \infty)}(0)$ is an isomorphism as a map $C^{2, \alpha}(L \oplus N_-) \lra C^\alpha(L \oplus N_-)$ for $\alpha \in (0, 1)$. Note that $N_{(\textbf{u}, \cdot)}$ is a smooth map from a neighborhood of $0$ in $C^{2, \alpha}(L \oplus N_-) \times \R$ to $C^\alpha (L \oplus N_-)$. The Banach space implicit function theorem gives the existence of a smooth map $f_{\frac{1}{R}}: [0, \frac{1}{R_0}) \lra C^{2, \alpha}(L \oplus N_-)$ such that $N_{(\textbf{u}, R)}(f_{\frac{1}{R}}) = 0$ for $0 \leq \frac{1}{R} < \frac{1}{R_0}$ and $f_0 = 0$. Undoing the change of variable we have a smooth map $f_R : (R_0, \infty) \lra C^{2, \alpha}(L \oplus N_-)$ such that $N_{(\textbf{u}, R)}(f_R) = 0$ for $R > R_0$ and $f_R \to 0$ as $R \to \infty$. Since all data in $N_{(\textbf{u}, R)}$ is real analytic, the real analytic version of the implicit function theorem \cite{MR2977424} shows that $f_R$ is real analytic in $R$. Finally, we need to check that $f_R \in \End_{h_R'}(E)$. Indeed,
    \begin{eqnarray*}
        (N_{(\textbf{u}, R)}(f_R))^{\ast_{h_R'}} = 
        e^{2 (f_R)^{\ast_{h_R'}}} N_{(\textbf{u}, R)}(f_R^{\ast_{h_R'}}) e^{-2 f_R^{\ast_{h_R'}}}.
    \end{eqnarray*}
    This shows that if $f_R$ is a solution, then $f_R^{\ast_{h_R'}}$ must be a solution as well. By the uniqueness part in the implicit function theorem we get $f_R = f_R^{\ast_{h_R'}}$, i.e., $f_R \in \End_{h_R'}(E)$, which completes the proof.
\end{proof}

\begin{proof} (of Theorem \ref{CL})
Let $[(\bp_E, \Phi)] \in \VHS_\alpha$ be a stable point. We first modify the harmonic metric $h$ to get a metric $h_R'$ by using (\ref{Gauge}). By Lemma \ref{LCL1}, we have
\begin{eqnarray*}
    \tilde{D}_{(\textbf{u}, R)} = \hbar^{-1} \Phi_\textbf{u} + \bp_\textbf{u} + \p_E - \dsp\sum_{j = -1}^{1 - \ell} R^{2 j} (\beta_j)^\ast + \hbar \Phi^\ast + \hbar \dsp\sum_{j = -2}^{1 - \ell} R^{2j + 2} (\varphi_j)^\ast.
\end{eqnarray*}
Taking the limit as $R \to \infty$ we obtain
\begin{eqnarray*}
    D_{(\textbf{u}, \infty)} := \dsp \lim_{R \to \infty} \tilde{D}_{(\textbf{u}, R)} = \hbar^{-1} \Phi_\textbf{u} + \bp_\textbf{u} + \p_E + \hbar \Phi^\ast.
\end{eqnarray*}
Observe that the connections $\tilde{D}_{(\textbf{u}, R)}$ are not flat, but they have the desired limit as $R \to \infty$.

By Lemma \ref{LCL2}, we can write the harmonic metric for $(\bp_\textbf{u}, R \Phi_\textbf{u})$ as $h(\textbf{u}, R) = g(h_R')$ with $g = \exp(f_R)$, for $f_R$ a (traceless) $h_R'$-hermitian endomorphism, for $R > 0$ large. Then $k = \exp(2 f_R)$ in equations (\ref{1}), (\ref{2}) and (\ref{3}), and we have
\begin{eqnarray*}
    \dsp \lim_{R \to \infty} D_{(\textbf{u}, R)} &=& \dsp \lim_{R \to \infty} \lb \hbar^{-1} \Phi_\textbf{u} + \bp_\textbf{u} + \p_\textbf{u}^{h_{(\textbf{u}, R)}} + \hbar R^2 \Phi_\textbf{u}^{\ast_{h_{(\textbf{u}, R)}}} \rb\\
    &=& \dsp \lim_{R \to \infty} \lb \hbar^{-1} \Phi_\textbf{u} + \bp_\textbf{u} + e^{-2 f_R} \circ \p_E^{h_R'} \circ e^{2 f_R} + \hbar e^{-2 f_R} (R^2 \Phi_\textbf{u}^{\ast_{h_R'}}) e^{2 f_R} \rb,
\end{eqnarray*}
which is equal to $D_{(\textbf{u}, \infty)}$ since $f_R \to 0$ as $R \to \infty$ by Lemma \ref{LCL3}.

\end{proof}

\subsection{Conformal limits in the Deligne-Hitchin moduli space}

Motivated by the following example in rank 2, we are interested in studying the realationship between the upward conformal limit and the conformal limit in the Deligne-Hitchin moduli space.
\begin{Example} \label{Ex}
    Consider the uniformizing Higgs bundle of rank two from (\ref{Fuchsian2}), which is a stable fixed point in $\Mm_\sH$. The Higgs bundles that flow up to it are of the form
    \begin{eqnarray} \label{E 4}
        \left( \begin{pmatrix}
        \bp_{K^{\frac{1}{2}}} & 0\\
        \beta & \bp_{K^{-\frac{1}{2}}}
        \end{pmatrix}, \begin{pmatrix}
        0 & 0\\
        1 & 0
        \end{pmatrix}\right),
    \end{eqnarray}
    with $\beta \in H^1(K^{-1})$ (this means that $\beta$ is harmonic: $\p_E \beta = 0$). In particular, all of them are in the nilpotent cone. By Theorem \ref{CL}, their upward conformal limit (we consider $\hbar = 1$ for simplicity) exists and is the flat connection
    \begin{eqnarray*}
        \left[ \left( 1,
        \begin{pmatrix}
            \bp_{K^{\frac{1}{2}}} & 0\\
            \beta & \bp_{K^{-\frac{1}{2}}}
        \end{pmatrix} + \begin{pmatrix}
            0 & 1^\ast\\
            0 & 0
        \end{pmatrix}, \begin{pmatrix}
            \p_{K^{\frac{1}{2}}} & 0\\
            0 & \p_{K^{-\frac{1}{2}}}
        \end{pmatrix} + \begin{pmatrix}
            0 & 0\\
            1 & 0
        \end{pmatrix} \right)\right].
    \end{eqnarray*}
    Recalling Definition \ref{DeH}, in the Deligne-Hitchin moduli space this 1-connection on $X$ is equal to the following 1-connection on $\Bar{X}$:
    \begin{eqnarray*}
        \left[ \left( 1, \begin{pmatrix}
            \p_{K^{-\frac{1}{2}}} & 0\\
            0 & \p_{K^{\frac{1}{2}}}
        \end{pmatrix} + \begin{pmatrix}
            0 & 1\\
            0 & 0
        \end{pmatrix},
        \begin{pmatrix}
            \bp_{K^{-\frac{1}{2}}} & 0\\
            0 & \bp_{K^{\frac{1}{2}}}
        \end{pmatrix} + \begin{pmatrix}
            0 & \beta\\
            1^\ast & 0
        \end{pmatrix} \right)\right].
    \end{eqnarray*}
    Notice that it is the image of $ \beta \in H^0(K^2)$ through the map $p_{\oper}^{\Bar{X}}$ from (\ref{oper}). Thus, the upward conformal limit of (\ref{E 4}) is an $\sSL(2, \C)$-oper in $\Bar{X}$. Hence, in rank two, the upward conformal limit of points that flow up to the zero of the Hitchin section on $X$ is an oper on $\Bar{X}$.
\end{Example}
This example holds in any rank, which gives Corollary \ref{MT3} (see Corollary \ref{C 4} below). First, we will give a more general identification between the filtrations of the nilpotent cone on $X$ and the Morse strata on $\Bar{X}$ described in Section \ref{DF}.

\begin{Theorem} \label{T 2}
    Let $\VHS_\alpha \subset \Mm_\sH^X$ be a connected component of the set of fixed points. For any $1 \leq j \leq \ell - 1$,
    \begin{eqnarray*}
        CL_{\hbar, \infty}^X(\NC_\alpha^{-j}) = CL_{\hbar^{-1},0}^{\Bar{X}}(W_\alpha^{0, j}).
    \end{eqnarray*}
\end{Theorem}

\begin{Remark} \label{R 3}
    Similarly, we can identify the filtration of a component of the Morse strata of $X$ with the filtration an irreducible component of the nilpotent cone on $\Bar{X}$:
    \begin{eqnarray*}
        CL_{\hbar, 0}^X(W_\alpha^{0, j}) = CL_{\hbar^{-1},\infty}^{\Bar{X}}(\NC_\alpha^{-j}).
    \end{eqnarray*}
\end{Remark}

The proof requires the following:

\begin{Lemma} \label{L 2}
    Let $[(\bp_0, \Phi_0)] \in \VHS_\alpha$ be stable and let $(\beta, \varphi) \in S^+_{(\bp_0, \Phi_0)}$. For any $R > 0$,
    \begin{eqnarray*}
        CL_{\hbar, 0}([\bp_0 + \beta, \Phi_0 + \varphi]) = CL_{\hbar R^{-1}, 0}([\bp_0 + \tilde{\beta}, \Phi_0 + \tilde{\varphi}])
    \end{eqnarray*}
    where $(\tilde{\beta}, \tilde{\varphi}) = R^{-1} \cdot (\beta, \varphi) \in S^+_{(\bp_0, \Phi_0)}$ is the action from (\ref{E 2}).
\end{Lemma}

\begin{proof}
    Let $R > 0$. Consider the gauge transformation $g_{R^{-1}}$ from (\ref{Gauge}). Then,
    \begin{eqnarray*}
        (1, \bp_0 + \beta + \hbar \Phi_0^\ast, \p_0 + \hbar^{-1}(\Phi + \varphi)) \cdot g_{R^{-1}} = (1, \bp_0 + \tilde{\beta} + \hbar R^{-1} \Phi_0^\ast, \p_0 + \hbar^{-1}R \Phi_0 + \hbar^{-1}R \tilde{\varphi}),
    \end{eqnarray*}
    where $(\tilde{\beta}, \tilde{\varphi}) := \dsp \left(\sum_{j=1}^{\ell -1} R^{-j} \beta_j, \sum_{j=0}^{\ell - 1} R^{-j-1} \varphi_j\right)$. In order to write the right hand side as the conformal limit of a Higgs bundle we need to check that $(\tilde{\beta}, \tilde{\varphi}) \in S^+_{(\bp_0, \Phi_0)}$. Indeed, observe that $(\tilde{\beta}, \tilde{\varphi}) = R^{-1} \cdot (\beta, \varphi)$ with $(\beta, \varphi) \in S^+_{(\bp_0, \Phi_0)}$, and by Lemma 3.10 in \cite{MR3949609}, the $\C^\ast$-action on $\Omega^{0,1}(\End(E)) \oplus \Omega^{1,0}(\End(E))$ from (\ref{E 2}) preserves the $\rB \rB$-slice. Thus, we can write
    \begin{eqnarray*}
        CL_{\hbar, 0}([\bp_E + \beta, \Phi + \varphi]) 
        = CL_{\hbar R^{-1}, 0}([\bp_0 + \tilde{\beta}, \Phi_0 + \tilde{\varphi}]).
    \end{eqnarray*}
\end{proof}

\begin{Remark}
    The same result as in Lemma \ref{L 2} holds true for the upward conformal limit of Higgs bundle that flows up to a stable $[(\bp_\infty, \Phi_\infty)] \in \VHS_\alpha$.
\end{Remark}

\begin{proof}(of Theorem \ref{T 2})
    Let $1 \leq j_0 \leq \ell - 1$ be fixed and let $[(\bp_E, \Phi)] \in \NC_\alpha^{-j_0}$ be stable. Recall that $\NC_\alpha$ can be characterized as $W_\alpha^\infty$ \cite{MR1650276, MR4519161}, and $\NC_\alpha^{-j_0} = W_\alpha^{\infty, -j_0}$. So there exists $[(\bp_\infty, \Phi_\infty)] \in \VHS_\alpha$ stable such that $\dsp \lim_{\xi \to \infty} \xi \cdot [(\bp_E, \Phi)] = [(\bp_\infty, \Phi_\infty)]$. In particular, $[(\bp_E, \Phi)] \in W_\alpha^{\infty, -j_0}(\bp_\infty, \Phi_\infty)$ and one can write $\bp_E = \bp_\infty + \beta$ and $\Phi = \Phi_\infty + \varphi$ for some $(\beta, \varphi) \in S^{\leq -j_0}_{(\bp_\infty, \Phi_\infty)}$, by Proposition \ref{P 1} and Corollary \ref{C 2}. Then the upward conformal limit of $[(\bp_E, \Phi)]$ exists (Theorem \ref{CL}) and is
    \begin{eqnarray*}
        CL_{\hbar, \infty}^X([(\bp_\infty + \beta, \Phi_\infty + \varphi)]) = [(1, \bp_\infty + \beta + \hbar \Phi_\infty^\ast, \p_\infty + \hbar^{-1}(\Phi_\infty + \varphi))],
    \end{eqnarray*}
    where $\p_\infty := \p_\infty^h$ and $\Phi_\infty^\ast$ are the (1,0)-part of the Chern connection and the adjoint with respect to the harmonic metric $h$ of $(\bp_\infty, \Phi_\infty)$ respectively. Using the gluing from the Deligne-Hitchin moduli space (Definition \ref{DeH}) we get a flat connection on $\Bar{X}$:
    \begin{eqnarray*}
        [(1, \bp_\infty + \beta + \hbar \Phi_\infty^\ast, \p_\infty + \hbar^{-1}(\Phi_\infty + \varphi))] = [(1, \p_\infty + \hbar^{-1}(\Phi_\infty + \varphi), \bp_\infty + \beta + \hbar \Phi_\infty^\ast)].
    \end{eqnarray*}
    Our goal is to find a Higgs bundle $p$ in $\Bar{X}$ and an element $\nu$ in its $\rB\rB$-slice so that $p + \nu \in \Ll_\alpha^{\Bar{X}}$ and the right hand side of the above expression is the conformal limit at zero of $p + \nu$.\\
    By Lemma \ref{L 4}, $(\p_\infty, \Phi_\infty^\ast)$ is a stable complex variation of Hodge structure on $\Bar{X}$ and
    \begin{eqnarray*}
        (\varphi, \beta) \in \Omega^{0,1}(N_{> j_0}^{\Bar{X}}) \oplus \Omega^{1,0}(N_{\geq j_0}^{\Bar{X}}).
    \end{eqnarray*}
    However, $(\p_\infty + \varphi, \Phi_\infty^\ast + \beta)$ is not a Higgs bundle on $\Bar{X}$ because $(\varphi, \beta) \not \in S_{(\p_\infty, \Phi_\infty^\ast)}^+$ since $[\varphi, \beta]$ may not be zero in general.
    \begin{itemize}
        \item Case 1: Assume that $\hbar \in U(1)$. The Higgs bundle $(\p_\infty, \hbar  \Phi_\infty^\ast)$ is stable and a complex variation of Hodge structure with harmonic metric $h$ and associated flat connection
        \begin{eqnarray*}
            D := \p_\infty + \bp_\infty + \hbar \Phi_\infty^\ast + \hbar^{-1} \Phi_\infty.
        \end{eqnarray*}
        Moreover, $D + \hbar^{-1} \varphi + \beta$ is flat with $(\hbar^{-1} \varphi, \beta) \in \Omega^{0,1}(N_{> j_0}^{\Bar{X}}) \oplus \Omega^{1,0}(N_{\geq j_0}^{\Bar{X}})$. By Propostion 4.8 in \cite{MR3949609} and Proposition \ref{P 3}, there is a unique $f \in \Omega^0(N_{>j_0}^{\Bar{X}})$ such that applying the complex gauge transformation $g = \mathbf{I} + f$,
        \begin{eqnarray*}
            (D + \hbar^{-1} \varphi + \beta) \cdot g = D + \tilde{\varphi} + \tilde{\beta},
        \end{eqnarray*}
        with $(\tilde{\varphi}, \tilde{\beta}) \in S^{\geq j_0}_{(\p_\infty, \hbar \Phi_\infty^\ast)}$. One checks that $(\tilde{\varphi}, \hbar^{-1} \tilde{\beta}) \in S^{\geq j_0}_{(\p_\infty, \Phi_{\infty}^\ast)}$ since $\hbar \in U(1)$. Thus,
        \begin{eqnarray*}
            CL_{\hbar, \infty}^X([(\bp_\infty + \beta, \Phi + \varphi)]) 
            &=& [(1, \p_\infty + \hbar^{-1}(\Phi_\infty + \varphi), \bp_\infty + \beta + \hbar \Phi_\infty^\ast)]\\
            &=& [(1, \p_\infty + \hbar^{-1}\Phi_\infty + \tilde{\varphi}, \bp_\infty + \hbar \hbar^{-1} \tilde{\beta} + \hbar \Phi_\infty^\ast)]\\
            &=& CL_{\hbar^{-1}, 0}^{\Bar{X}}([(\p_\infty + \tilde{\varphi}, \Phi_\infty^\ast + \hbar^{-1} \tilde{\beta})]),
        \end{eqnarray*}
        where the last limit makes sense. Moreover, Lemma \ref{L 3} guarantees that on $\Bar{X}$
        \begin{eqnarray*}
            (\p_\infty + \tilde{\varphi}, \Phi_\infty^\ast + \hbar^{-1} \tilde{\beta}) \in W_\alpha^{0, j_0}(\bp_\infty, \Phi_\infty).
        \end{eqnarray*}

        \item Case 2: Suppose now $\hbar \not \in U(1)$. We can write $\hbar = |\hbar| e^{i \theta}$ for some $\theta \in [0, 2\pi]$. Applying the gauge transformation $g_{|\hbar|}$ from (\ref{Gauge}) to $(\p_\infty, \hbar \Phi_\infty^\ast)$ we get the Higgs bundle $(\p_\infty, e^{i \theta} \Phi_\infty^\ast)$. By Case 1, there is $(\tilde{\varphi}, e^{-i \theta} \tilde{\beta}) \in S^{\geq j_0}_{(\p_\infty, \Phi_\infty^\ast)}$ such that
        \begin{eqnarray*}
            CL_{\hbar, \infty}^X([(\bp_\infty + \beta, \Phi + \varphi)]) 
            &=& [(1, \p_\infty + \hbar^{-1}(\Phi_\infty + \varphi), \bp_\infty + \beta + \hbar \Phi_\infty^\ast)]\\
            &=& [(1, \p_\infty + e^{-i \theta} \Phi_\infty + \tilde{\varphi}, \bp_\infty + \tilde{\beta} + e^{i \theta} \Phi_\infty^\ast)]\\
            &=& CL^{\Bar{X}}_{e^{-i \theta}, 0}([(\p_\infty + \tilde{\varphi}, \Phi_\infty^\ast + e^{-i \theta} \tilde{\beta})]).
        \end{eqnarray*}
        Undoing the gauge transformation to recover the original Higgs bundle $(\p_\infty, \hbar \Phi_\infty^\ast)$ and using Lemma \ref{L 2} we get
        \begin{eqnarray*}
            CL^{\Bar{X}}_{e^{-i \theta}, 0}([(\p_\infty + \tilde{\varphi}, \Phi_\infty^\ast + e^{-i \theta} \tilde{\beta})]) = CL^{\Bar{X}}_{\hbar^{-1}, 0}([(\p_\infty + \varphi, \Phi_\infty^\ast + e^{-i \theta} \beta)]),
        \end{eqnarray*}
        where $(\varphi, e^{-i\theta} \beta) := |\hbar|^{-1} \cdot (\tilde{\varphi}, e^{-i \theta} \tilde{\beta}) \in S^{\geq j_0}_{(\p_\infty, \Phi_\infty^\ast)}$. Note that $(\p_\infty + \varphi, \Phi_\infty^\ast + e^{-i \theta} \beta)$ is an element of $W_\alpha^{0, j_0}(\bp_\infty, \Phi_\infty)$ on $\Bar{X}$ by Lemma \ref{L 3}.
    \end{itemize}
\end{proof}

A family of Lagrangians on $\Mm_\sH$ was recently introduced in \cite{gothen2026zerovelocitylagrangianshitchin}. They are called \textit{zero velocity Lagrangians} because they consist of points which ``arrive at their limit to 0 with zero velocity". Namely,
\begin{eqnarray*}
    \Ll_\alpha := \{ [(\bp_E, \Phi)] \in W^0_\alpha : \gamma'_{[(\bp_E, \Phi)]} (0) = 0\},
\end{eqnarray*}
where
\begin{eqnarray*}
    \begin{array}{rc c l}
      \gamma_{[(\bp_E, \Phi)]} :& \C & \longrightarrow & \Mm_\sH  \\
      &0 \neq t & \longmapsto & t \cdot [(\bp_E, \Phi)]\\
      &0 & \longmapsto & \dsp \lim_{t \to 0}t \cdot[(\bp_E, \Phi)].
    \end{array} 
\end{eqnarray*}
The Higgs bundles on $\Ll_\alpha$ are called zero velocity Higgs bundles and can be characterized using the filtration of the $\rB \rB$-slice as in the following result \cite[Proposition 3.5]{gothen2026zerovelocitylagrangianshitchin}.

\begin{Proposition} \label{P 4}
    Let $[(\bp_E, \Phi)] \in \Mm_\sH$ such that $\dsp \lim_{t \to 0}t \cdot[(\bp_E, \Phi)] = [(\bp_0, \Phi_0)] \in \VHS_\alpha$ is stable. Then $[(\bp_E, \Phi)]$ is a zero velocity Higgs bundle if and only if it can be represented as $[(\bp_0 + \beta, \Phi_0 + \varphi)]$ for some pair $(\beta, \varphi) \in S^{\geq 1}_{(\bp_0, \Phi_0)}$.
\end{Proposition}

By Lemma \ref{L 3}, we can write $\Ll_\alpha = W_\alpha^{0, 1}$. As a consequence of Theorem \ref{T 2}, in the Deligne-Hitchin moduli space a connection between zero velocity Lagrangians on $\Bar{X}$ and irreducible components of the nilpotent cone on $X$ arises when we look at their images under the conformal limit and the upward conformal limit, respectively. This gives us Theorem \ref{MT2}:

\begin{Corollary}
    Let $\VHS_\alpha \subset \Mm_\sH^X$ be a connected component of the set of fixed points. Let $\NC_{\alpha}^X$ be the irreducible component of the nilpotent cone on $X$ and $\Ll_{\alpha}^{\Bar{X}}$ be the zero velocity Lagrangian submanifold of the Morse strata on $\Bar{X}$, both associated to $\VHS_\alpha$. Then,
    \begin{eqnarray*}
        \CL_{\hbar, \infty}^{X}(\NC_{\alpha}) = \CL_{\hbar^{-1}, 0}^{\Bar{X}}(\Ll_{\alpha}).
    \end{eqnarray*}
\end{Corollary}

Now we can generalize Example \ref{Ex}. Observe that in the rank two case, the Higgs bundles of the form (\ref{E 4}) are all the points that flow up to the Fuchsian one. However, in higher rank we need to be more careful and we cannot consider all the points in the unstable manifold. In general, if we consider all the points that flow up to any stable fixed point, we obtain partial opers on $\Bar{X}$. 

\begin{Corollary} \label{C 4}
    Let $[(\bp_\infty, \Phi_\infty)]$ be the $\C^\ast$-fixed point in a Hitchin section on $X$. Then
    \begin{eqnarray*}
    \NC_\alpha(\bp_\infty, \Phi_\infty) \cong \bigoplus_{j = 1}^{n - 1} H^0(K_{\Bar{X}}^{j + 1}).
    \end{eqnarray*}
    Moreover, for any $\beta \in \dsp \bigoplus_{j = 1}^{n-1} H^1(K_{X}^{-j})$, the $\hbar$-upward conformal limit of the Higgs bundle associated to $\beta$ is the oper on $\Bar{X}$ associated to $\beta$ viewed as an element of $\dsp \bigoplus_{j = 1}^{n-1} H^0(K_{\Bar{X}}^{j + 1})$.
\end{Corollary}
    
\begin{proof}
    Recall that $W_\alpha^0(\bp_\infty, \Phi_\infty)$ is parametrized by the vector space $\dsp \bigoplus_{j = 1}^{n - 1} H^0(K_{X}^{j + 1})$. By Remark \ref{R 3}, we have that $W_\alpha^\infty(\bp_\infty, \Phi_\infty)$ is parametrized by
    \begin{eqnarray*}
        \bigoplus_{j = 1}^{n-1} H^1(K_{X}^{-j}) \cong \bigoplus_{j = 1}^{n - 1} H^0(K_{\Bar{X}}^{j + 1}).
    \end{eqnarray*}
    
    Let $\beta \in \dsp \bigoplus_{j = 1}^{n-1} H^1(K_{X}^{-j})$. The Higgs bundle $[(\bp_\infty + \beta, \Phi_\infty)]$ flows up to the $\C^\ast$-fixed point. Then, by direct computation we have
    \begin{eqnarray*}
        \CL^X_{\hbar,\infty}([(\bp_\infty + \beta, \Phi_\infty)]) = \CL^{\Bar{X}}_{\hbar^{-1},0}([(\p_\infty, \Phi_\infty^\ast + \beta)]),
    \end{eqnarray*}
    since in this case $(\p_\infty, \Phi_\infty^\ast + \beta)$ is a Higgs bundle on $\Bar{X}$. Note that $\beta$ is a quadratic differential on $\Bar{X}$, and the right hand side of the above equality is the image of $\beta$ through the map $p_{\oper}$ from (\ref{oper}). Hence, $\CL^X_{\hbar,\infty}([(\bp_\infty + \beta, \Phi_\infty)])$ is an oper on $\Bar{X}$.
\end{proof}

\bibliographystyle{plain}
\bibliography{References}{}

\end{document}